\documentclass[a4paper,11pt]{article}
\usepackage[nointlimits]{amsmath}
\usepackage{amsfonts}
\usepackage{amssymb}
\usepackage{tensor}
\usepackage{mathtools}
\usepackage{amsfonts}
\usepackage{amssymb}
\usepackage{caption}
\usepackage{subcaption}
\usepackage{tikz}
\usepackage{amsthm}
\usepackage{amsmath}
\usepackage{latexsym}
\usepackage{graphicx}
\usepackage{layout}
\usepackage[english]{babel}
\numberwithin{equation}{section}

\usepackage{graphics}
\usepackage{epsfig}

\usepackage{tikz}
\usepackage{color}
\usepackage{graphicx}
\usepackage{caption}
\usepackage{subcaption}

\newtheorem{lemma1}     {Lemma}[section]
\newtheorem{teorema1}   [lemma1]{Theorem}
\newtheorem{prop1}      [lemma1]{Proposition}
\newtheorem{coroll1}    [lemma1]{Corollary}
\newtheorem{cong1}      [lemma1]{Conjecture}
\newtheorem{remark1}    [lemma1]{Remark}
\newtheorem{defin1}     [lemma1]{Definition} 
\newtheorem{example1}   [lemma1]{Example} 

\newenvironment{Lemma}[1][]
        {\begin{lemma1}[#1]\begin{samepage}}{\end{samepage}\end{lemma1}}
\newenvironment{Theorem}[1][]
        {\begin{teorema1}[#1]\begin{samepage}}{\end{samepage}\end{teorema1}}
\newenvironment{Proposition}[1][]
        {\begin{prop1}[#1]\begin{samepage}}{\end{samepage}\end{prop1}}
\newenvironment{Corollary}[1][]
        {\begin{coroll1}[#1]\begin{samepage}}{\end{samepage}\end{coroll1}}

\newenvironment{Remark}[1][]
        {\begin{remark1}[#1]\begin{samepage}}{\end{samepage}\end{remark1}}
\newenvironment{Definition}[1][]
        {\begin{defin1}[#1]\begin{samepage}}{\end{samepage}\end{defin1}}

\newcommand{\nada}[1]   {}
\newcommand{\A}         {\mathcal A}
\newcommand{\Aa}{{E}_1}
\newcommand{\Aae}{{E}_1^\eps}
\newcommand{\ac}{{23}}
\newcommand{\alp}{{\alpha}_1}
\newcommand{\area}         {\mathbb A}

\newcommand{\Bb}{{E}_2}
\newcommand{\Bbe}{{E}_2^\eps}
\newcommand{\bc}{{31}}
\newcommand{\bet}{{\alpha}_2}
\newcommand{\C}         {\mathcal C}
\newcommand{\calpbet}{C_{12}}

\newcommand{\cc}{{12}}
\newcommand{\Cc}{{E}_3}
\newcommand{\Cce}{{E}_3^\eps}

\newcommand{\cij}{C_{ij}}

\newcommand{\currentijn}{{[\![ \pnGammaalpij ]\!]}}
\newcommand{\currentijb}{{[\![ \bGammaalpij ]\!]}}
\newcommand{\currentij}{{[\![ \Gammaalpij ]\!]}}

\newcommand{\currentki}{{[\![ \Gammaalpki ]\!]}}

\newcommand{\currentib}{{[\![ \bGammaalpi ]\!]}}

\newcommand{\DD}	   {\mathcal D}

\newcommand{\del}{\delta}
\newcommand{\dirdatum}      {\varphi}

\newcommand{\eps}       {\varepsilon}
\newcommand{\eti}{\zetai}
\newcommand{\etii}{\zetaii}
\newcommand{\etiii}{\zetaiii}
\newcommand{\etia}{\zetai_1}
\newcommand{\etiia}{\zetaii_1}

\newcommand{\G}         {\mathcal G}
\newcommand{\gam}{{\alpha}_3}
\newcommand{\Gammacc}{\calpbet}
\newcommand{\Gammaalpij}{\tensor[]{\Gamma}{_{{ij}}}}

\newcommand{\Gammaalpki}{\tensor[]{\Gamma}{_{\ki}}}

\newcommand{\genericpn}{{\genericp_n}}

\newcommand{\grad}      {\nabla}
\renewcommand{\H}       {\mathcal{H}}

\newcommand{\Hone}        {\ensuremath{\mathcal H^1}}
\newcommand{\Lip}{{\rm Lip}} 

\newcommand{\N}         {\ensuremath{\mathbb N}}

\newcommand{\phiij}{\dirdatum_{i j}}

\newcommand{\raggio}       {r}

\newcommand{\rettangolo}       {\mathrm R}
\newcommand{\rectangle}{\rettangolo}

\newcommand{\R}         {\ensuremath{\mathbb R}}

\newcommand{\scalpbet}{c_{12}}
\newcommand{\scij}{c_{ij}}
 
\newcommand{\tenda}       {m} 
 
\newcommand{\Triangle}{\mathrm T}
\newcommand{\uh}       {u^\eps} 
\newcommand{\uhs}       {u^{\eps}}
\newcommand{\ww}{w}
\newcommand{\zetai}{\zeta^1}
\newcommand{\zetaii}{\zeta^2}
\newcommand{\zetaiii}{\zeta^3}

\newcommand{\sGammacc}{\scalpbet}

\newcommand{\e}{e}
\newcommand{\g}{g}
\newcommand{\NN}{\mathcal{N}}
\newcommand{\NNiii}{\mathcal{N}_\cc^{\eps,\del}}

\newcommand{\ro}{\rho}

\newcommand{\h}{h}
\newcommand{\z}{\tau}
\newcommand{\w}{d}
\newcommand{\nub}{\bar{\nu}}
\newcommand{\hinv}{\h^{-1}}
\newcommand{\neworigin}{Q}
\newcommand{\ce}{c_\eps}

\newcommand{\Xc}{X}
\newcommand{\Xclip}{X_{\Lip}}

\newcommand{\elementXclip}{ \Gamma}
\newcommand{\barelementXclip}{ \widetilde\Gamma}

\newcommand{\Lone}{L^1} 

\newcommand{\sourcedisk}{D}
\newcommand{\targettriangle}{\mathrm T}
\newcommand{\sourceorigin}{0_S}
\newcommand{\targetorigin}{0_{\targettriangle}}

\newcommand{\genericp}{p}

\newcommand{\alpi}{\alpha_i}
\newcommand{\alpj}{\alpha_j}

\newcommand{\ki}{ki}
\newcommand{\pGammaalpij}{\Gamma_{{ij}}}

\newcommand{\bGammaalpij}{\Gamma_{{ij}}^n}

\newcommand{\pnGammaalpij}{\Gamma_{{ij}}^n}
\newcommand{\pnGammaalpki}{\Gamma_{\ki}^n}

\newcommand{\bGammaalpi}{\Gamma_{i}^n}
\newcommand{\pbarGammaalpi}{\widetilde\Gamma_{i}}

\newcommand{\Gammaalp}{\Gamma_{1}}
\newcommand{\Gammabet}{\Gamma_{2}}
\newcommand{\Gammagam}{\Gamma_{3}}
\newcommand{\pGammaalpi}{\Gamma_{i}}
\newcommand{\pGammaalpj}{\Gamma_{j}}

\newcommand{\ellalpbet}{\ell_{1 2}}
\newcommand{\ellbetgam}{\ell_{2 3}}
\newcommand{\ellgamalp}{\ell_{3 1}}

\newcommand{\phialpbet}{\dirdatum_{1 2}}

\newcommand{\widehatphialpbet}{\widehat{ \dirdatum}}

\newcommand{\phialpij}{\dirdatum_{i j}}
\newcommand{\phialpki}{\dirdatum_{k i}}

\newcommand{\wo}{w}
\newcommand{\wij}{w_{i j}}
\newcommand{\wki}{w_{k i}}

\newcommand{\ellalpij}{{\ell_{i j}}}
\newcommand{\ellalpki}{\ell_{k i}}

\newcommand{\radiusalpij}{r_{i j}}

\newcommand{\radiusalpbet}{r_{12}}
\newcommand{\radiusbetgam}{r_{23}}
\newcommand{\radiusgamalp}{r_{31}}

\newcommand{\walpbet}{w_{1 2}}

\newcommand{\wgamalp}{w_{3 1}}

\newcommand{\Ball}{B}
\newcommand{\T}{\mathcal{T}}

\newcommand{\Tli}{\tensor[]{\T}{_i}}
\newcommand{\Tlj}{\tensor[]{\T}{_j}}
\newcommand{\Tlk}{\tensor[]{\T}{_k}}

\newcommand{\Ri}{\mathcal{P}_i}
\newcommand{\Rj}{\mathcal{P}_j}
\newcommand{\Rk}{\mathcal{P}_k}
\newcommand{\Ralp}{\mathcal{P}_1}

\newcommand{\Lij}{h_{ij}}
\newcommand{\Ljk}{h_{jk}}
\newcommand{\Lki}{h_{ki}}
\newcommand{\Lalpbet}{h_{12}}
\newcommand{\Lgamalp}{h_{31}}

\newcommand{\spt}{{~\rm spt~}}

\newcommand{\I}{(a,b)}
\newcommand{\cI}{[a,b]}
\newcommand{\BV}{{\rm BV}}

\newcommand{\AG}{\mathcal G}

\newcommand{\fe}{\varphi^n}
\newcommand{\ngamalp}{{n_{31}}}\newcommand{\nalpbet}{{n_{12}}}

\newcommand{\nualp}{\nu}

\newcommand{\keps}{\overline{\kappa}_\eps}

\newcommand{\hkindex}{h}

\newcommand{\res}{\mathop{\hbox{\vrule height 7pt width 0.5pt depth 0pt
\vrule height 0.5pt width 6pt depth 0pt}}\nolimits}

\graphicspath{{figures/}}

\begin{document}

\title{On the relaxed area of the graph of discontinuous maps from the plane to the plane taking three values with no symmetry assumptions 
}
\author{
Giovanni Bellettini\footnote{
Dipartimento di Ingegneria dell'Informazione e Scienze Matematiche, Universit\`a di Siena, 53100 Siena, Italy,
and International Centre for Theoretical Physics ICTP,
Mathematics Section, 34151 Trieste, Italy.
E-mail: bellettini@diism.unisi.it
                      }\and
Alaa Elshorbagy\footnote{
Area of Mathematical Analysis, Modelling, and Applications,
             Scuola Internazionale Superiore di Studi Avanzati "SISSA",
Via Bonomea, 265 - 34136 Trieste, Italy,
and 
International Centre for Theoretical Physics ICTP,
Mathematics Section, 34151 Trieste, Italy E-mail: alaa.elshorbagy@sissa.it
                         }
                          \and
Maurizio Paolini\footnote{
Dipartimento di Matematica e Fisica, Universit\`a Cattolica del Sacro Cuore, 25121 Brescia, Italy.
E-mail: paolini@dmf.unicatt.it
                         }
 \and
Riccardo Scala\footnote{ Dipartimento di Matematica ``Guido Castelnuovo'', Universit\`a La Sapienza, Piazzale Aldo Moro 5, 00185 Roma.
E-mail: scala@mat.uniroma1.it }
}

\date{}
\maketitle
\thanks{}

\begin{abstract}
In this paper we estimate  from above the area of the
graph of a  singular map $u$ taking a disk to  three
vectors,  the vertices of a triangle,  and
jumping along three $\C^2-$ embedded curves that meet transversely at only
one point of the disk.
We show that  the relaxed area  can be estimated from above by the
solution of a Plateau-type problem involving three entangled
nonparametric area-minimizing surfaces.
The idea is to ``fill
the hole'' in the graph of the singular map
 with a sequence of approximating
smooth two-codimensional surfaces of
graph-type, by imagining three minimal surfaces, placed vertically over
the jump of $u$, coupled together via a triple point in the target
triangle. Such a construction depends on the choice of a target triple
point, and on a connection passing through it, which  dictate the
boundary condition for the three minimal surfaces.
We show that the singular part of the relaxed area of $u$
cannot be larger than what we obtain by minimizing over all possible
target triple points and all corresponding connections.
\end{abstract}

\noindent {\bf Key words:}~~Relaxation, Cartesian currents, area functional, minimal surfaces, Plateau problem.

\vspace{2mm}

\noindent {\bf AMS (MOS) subject clas\-si\-fi\-ca\-tion:} 
49Q15, 49Q20, 49J45.


\section{Introduction}\label{sec:introduction}
Let $\Omega \subset \R^2$ be an open set  and $v=(v_1,v_2) : \Omega \to  \R^2$
a Lipschitz map. It is well known that the area of 
the graph of $v$ is given by 
\begin{equation}\label{eq:area_functional}
 \area (v,\Omega)= \int_{\Omega}\sqrt{ 1+|\grad v_1|^2+|\grad v_2|^2+\Big( 
\frac{\partial v_1}{\partial x}\frac{\partial v_2}{\partial y}-\frac{\partial v_1}{\partial y}
\frac{\partial v_2}{\partial x}\Big)^2}dxdy. 
 \end{equation}
Extending to nonsmooth maps 
via relaxation the definition of the area
is a difficult question
\cite{Giq},
and is motivated by rather natural problems in calculus of variations: 
we can mention for example
the use of direct methods to face
the two-codimensional Plateau problem in $\R^4$ in cartesian form, 
and the study of lower semicontinuous
envelopes of polyconvex functionals with nonstandard growth \cite{AcDa:94}, 
\cite{FuHu:95}.
A crucial issue is to decide
which topology one has to consider in order to compute the relaxed 
functional of $\area(\cdot,\Omega)$: of course, the weakest the topology,
the most difficult should be the computation of the relaxed functional, but the easiest 
becomes the coerciveness.
 We recall 
that when $v$ is scalar valued, the natural choice is the $L^1(\Omega)$-convergence,
and 
the relaxation  problem is completely solved \cite{Da:80},
\cite{AmFuPa:00}; 
the $L^1(\Omega)$-relaxed functional in this case
consists, besides the absolutely continuous part, of a singular part which 
is  the total variation of the jump and Cantor
parts of the distributional derivative of $v$ in $\Omega$; in particular, 
the relaxed functional, when considered as a function of $\Omega$,
is a measure. 

The case of interest here, namely when $v$ takes values
in $\R^2$, is much more involved,
due to the {\it nonconvexity} of the integrand in \eqref{eq:area_functional},
and to the {\it unilateral linear growth} 
$$
\area (v,\Omega) \geq \int_\Omega \sqrt{\vert \nabla v_1\vert^2 + 
\vert \nabla v_2\vert^2} ~dxdy.
$$
Choosing again the 
$L^1(\Omega; \R^2)$-convergence
(as we shall do in this paper), the 
relaxed functional $\A(\cdot,\Omega)$ of $\area(\cdot,\Omega)$,
{\it i.e.},
\begin{equation}
\A(v, \Omega):=\inf \big \{ \liminf_{\epsilon \to 0}
\area (u^\epsilon,\Omega): \{u^\epsilon\} \subset 
{\rm Lip}(\Omega; \R^2), ~ u^\epsilon \to u~ {\rm in} ~L^1(\Omega; \R^2)\big \},
\end{equation}
is, for $v \in L^1(\Omega; \R^2) \setminus W^{1,2}(\Omega; \R^2)$, 
far from being understood, and  exhibits surprising features.
One of the few known facts that 
must be pointed out is that, for a large 
class of nonsmooth maps $v$, the function
$\Omega \to \A(v, \Omega)$ 
cannot be written as an integral \cite{AcDa:94}, \cite{BePaTe:15}, \cite{BePaTe:16};
this interesting phenomenon, related to 
nonlocality, has at least two sources. 
For simplicity, let us focus our attention on nonsmooth functions 
with jumps, thus neglecting 
the case of vortices.
The first source of nonlocality
has been enlightened answering to 
a conjecture in 
\cite{De:92}. Specifically, consider the  symmetric triple junction map $u_{\rm symm}$,
{\it i.e.}, the singular map from a disk $D$ of $\R_S^2=\R^2$ into $\R_T^2=\R^2$,
taking only three values -- the vertices of an equilateral triangle $T_{\rm eq}\subset \R^2_T$ -- and
jumping along three segments meeting 
at the origin in a triple junction
at equal $120^\circ$ angles: 
 then  
$\A(u_{\rm symm}, \cdot)$ is not subadditive. This result has been proven in 
\cite{AcDa:94}; subsequently 
in \cite{BePa:10} it is shown that the value
 $\A(u_{\rm symm}, D)$ 
is related to the solution of three one-codimensional 
Plateau-type problems in cartesian form suitably entangled together through the Steiner 
point in the triangle $T_{\rm eq}$. Due to the
special symmetry of the map $u_{\rm symm}$,  the three-problems
collapse together to only one one-codimensional Plateau-type problem in cartesian form,
on a fixed rectangle $R$ whose sides are the radius of $D$ and the side of
$T_{\rm eq}$.
 Positioning  three copies of 
this minimal surface ``vertically'' (in the space of graphs, {\it i.e.},
in $D \times \R^2$)  over the jump of $u_{\rm symm}$
allows, in turn, to construct
a sequence $\{u_\eps\}$ of {\it Lipschitz} maps from $D$ 
into $\R^2$ 
the limit area of which improves the upper estimate of \cite{AcDa:94}.
Optimality of this construction 
has been shown in 
the recent paper 
\cite{Sc:19}, on the basis of a symmetrization procedure for currents.

It is one of the aims of the present paper
to inspect solutions of the above mentioned three Plateau-type problems
in more general situations,
in order to provide upper estimates for 
$\A(u,D)$, for suitable piecewise constant maps $u$.

A second source of nonlocality for the functional $\A(u,\Omega)$ is given by the
interaction of the jump set of a discontinuous map $u$ with 
the boundary of the domain $\Omega$.
This phenomenon, already observed in \cite{AcDa:94} for the map
with one-vortex
at the center of  a suitable disk, appears also for functions with jump
discontinuities not piecewise constant \cite{BePaTe:16}. 
More surprisingly, it appears also for piecewise constant
maps taking three values, provided the jump is sufficiently
close to the boundary of $\Omega$, as observed in 
\cite{Sc:19}, taking as $\Omega$ a sufficiently thin 
tubular neighbourhood of the jump itself. We shall not be concerned here with this second
source of nonlocality. 

\begin{figure}
\centering
\begin{subfigure}{.5\textwidth }
  \centering
  \begin{center}
     \includegraphics[scale=0.4]{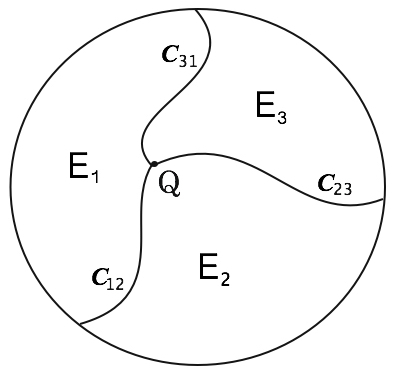}
   \end{center}
 \caption{The domain of $u$;
  $u=\alpha_i$ on $E_i$. }     \label{fig:sourcedisk}
  \end{subfigure}%
\begin{subfigure}{.5\textwidth} 
  \centering
  \begin{center}
       \includegraphics[scale=0.4]{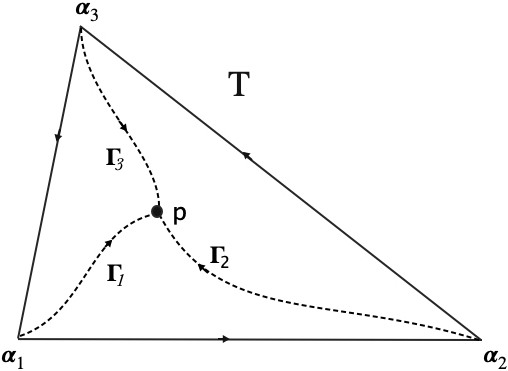}
  \end{center}
  \caption{A Lipschitz graph-type connection in the target triangle $\targettriangle$. 
$\Gamma_1 \cup \Gamma_2$ (resp. $\Gamma_2 \cup \Gamma_3$, $\Gamma_3 \cup \Gamma_1$)
is graph over the segment
$\overline{\alpha_1 \alpha_2}$
(resp.
$\overline{\alpha_2 \alpha_3}$,
$\overline{\alpha_3 \alpha_1}$) of a Lipschitz function 
$\varphi_{12}$
(resp. 
$\varphi_{23}$,
$\varphi_{31}$).
}\label{fig:targettriangle}
\end{subfigure}
\caption{}\label{fig:sourcedisk_targetriangle}
\end{figure}

As already mentioned above, 
in this paper we are interested in estimating from above
the area of the graph of a  singular map $u$ taking three 
(non collinear) values 
and jumping along three embedded curves of class $\C^2$ that meet transversely at only one point, see Figure \ref{fig:sourcedisk_targetriangle}.
Let us state this in a more precise way, referring to Sections
\ref{sec:notation} and \ref{sec:Connection} for all details.
For simplicity, from now on we fix $\Omega$ to be an open disk 
$\sourcedisk $ containing the origin $\sourceorigin$ in the source 
plane $\R^2= \R^2_{x,y}=\R_S^2$.
Take three non-overlapping non-empty two-dimensional connected regions $ \Aa,~\Bb,~\Cc$ 
of  $\sourcedisk $ 
such that 
\begin{equation}\label{eq:connected_regions}
  \Aa \cup \Bb \cup \Cc= \sourcedisk.
\end{equation}
The three regions are separated by three embedded curves of class $\C^2$ (up to the boundary) of length $\radiusalpbet, ~\radiusbetgam,~\radiusgamalp$ respectively, that meet only at $\neworigin$ (source triple junction); 
moreover, each curve is supposed to meet the 
boundary of $\sourcedisk$ transversely and we assume also that $\neworigin$ is a
transversal intersection for the 
three curves, see Figure \ref{fig:sourcedisk}.
Let $\alp,~ \bet, ~\gam$ be the vertices of a closed triangle $\targettriangle$ with non empty interior in the target plane.

Set
\begin{equation}\label{eq:l_ij}
\ellalpbet:=|\alp-\bet|, \qquad \ellbetgam:=|\bet-\gam|, \qquad \ellgamalp:=|\alp-\gam|.
\end{equation}
We suppose that 
$\targettriangle$ contains the origin $\targetorigin$ in its interior.

\begin{figure}
\centering
\begin{subfigure}{.5\textwidth }
  \centering
  \begin{center}
     \includegraphics[scale=0.34]{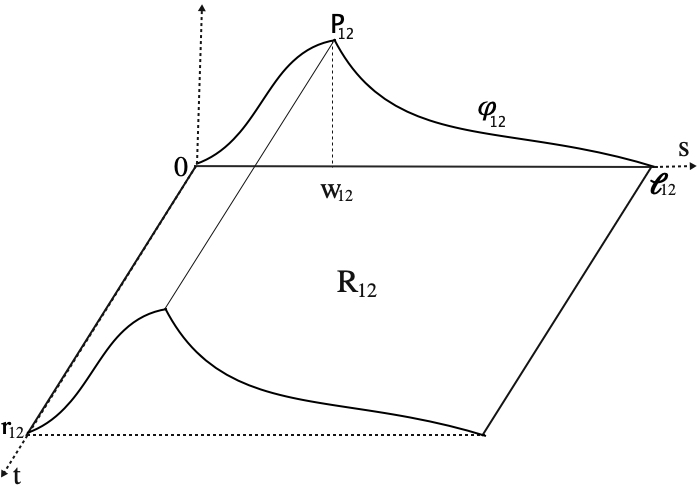}
   \end{center}
 \caption{The graph of the function $\dirdatum_\cc$ on $\rettangolo_{12}$.}  \label{fig:dirchlet_datum}
  \end{subfigure}%
\begin{subfigure}{.5\textwidth} 
  \centering
  \begin{center}
       \includegraphics[scale=0.34]{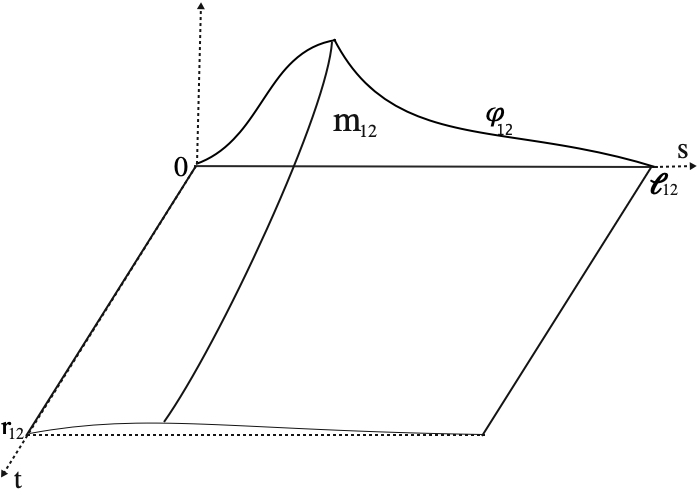}
  \end{center}
  \caption{The graph of $\tenda_\cc$ on $\rettangolo_{12}$.}  \label{fig:graph_m_i}
  \end{subfigure}
\caption{}\label{fig:rect_mini}
\end{figure}

Let us introduce the space $X$ of connections (Definition \ref{def:connections} and \eqref{eqn:classXlip}, \eqref{eqn:classX}); a connection 
$\Gamma = (\Gamma_1, \Gamma_2, \Gamma_3)$ consists of three rectifiable 
curves in $\targettriangle$,
that connect the vertices of $\targettriangle$
to some point 
inside the triangle (called target triple point). We shall suppose that {\it each curve
can be written as a graph, possibly with vertical parts, over
the corresponding two sides of $\targettriangle$}.  When $\Gamma$ 
consists of three Lipschitz graphs, we write $\Gamma \in X_{\rm Lip}$,
and we say that $\Gamma$ is a Lipschitz connection.
We  now show how to construct a new functional 
$\mathcal G$, consisting of 
the sum of the areas of three minimal surfaces -- graphs
of three suitable area-minimizing 
functions $\tenda_{\cc}$, $\tenda_{ \ac}$, $\tenda_{\bc}$ defined on certain
 rectangles -- coupled together 
by the connection considered as a Dirichlet 
boundary condition, see Definition \ref{def:functional_G}. Set 
\begin{equation}\label{eq:R_ij}
\rettangolo_{ij} := 
\left[0, \ellalpij \right] \times [0, r_{ij}],\qquad ij \in\{12,23,31\}.
\end{equation}
Assume $\elementXclip \in X$. Then 
$\pGammaalpij := \Gamma_i \cup \Gamma_j$, $ij\in\{12,23,31\} $ are (generalized) graphs of functions $\phiij$ of bounded variation over  $[0,\ellalpij]$. 
With a small abuse of notation,
set
 \begin{equation}\label{eq:extension_of_varphi}
 \phialpij (s,t)=\phialpij(s), \qquad (s,t)\in 
\rettangolo_{i j},  \qquad ij\in\{12,23,31\}.
 \end{equation}
The graph of $\phialpbet$ on $\rettangolo_{1 2}$ is depicted in Figure \ref{fig:dirchlet_datum}.

Let $\tenda_{ij}=\tenda_{i j}(\Gamma) $ be 
the unique solution of the Dirichlet-Neumann minimum problem
{
\begin{align}\label{eq:Dirichlet_Neumann_minimum_problem}
\min 
\left\{
\int_{\rettangolo_{i j}}
\sqrt{1 + \vert \grad f\vert^2} ~ds dt : ~ f \in W^{1,1}(\rettangolo_{i j}), 
~ f = \phialpij ~\H^1 -a.e. ~{\rm on}~ 
\partial_D \rettangolo_{i j}
\right\},
\end{align}
}
where
\begin{align*}
\partial_D \rettangolo_{i j}= \partial \rettangolo_{i j} \setminus \left( [0,\ellalpij] \times \{\radiusalpij\} \right), \qquad  ij\in\{12,23,31\} .
\end{align*}
Notice that 
the minimization is taken among all functions having a Dirichlet condition on three of the four sides of the rectangle $\rettangolo_{i j}$; 
the missing side corresponds to the intersection points of the jump
with the boundary of $D$.

From \eqref{eq:extension_of_varphi} it follows that the Dirichlet condition is zero on the sides $\{0\}\times [0,\radiusalpij]$ and $ \{\ellalpij\}  \times [0,\radiusalpij] $ of $\rettangolo_{i j}$;
see Figure \ref{fig:graph_m_i}.

Set 
\begin{equation}\label{eq:A_ij}
{\mathfrak A}_{i j} (\Gamma) 
:=  
\int_{\rettangolo_{i j}} \sqrt{1 + \vert\grad\tenda_{i j}\vert^2}~dsdt,\qquad  ij\in\{12,23,31\}.
\end{equation}
The main result of the present paper reads as follows (see Theorem 
\ref{thm:upper_bound_with_infimum} and Corollary \ref{COR}). 
\begin{Theorem}\label{teo:main}
Let $u: \sourcedisk \to \{\alp, \bet, \gam\}$ be the discontinuous
$BV(\sourcedisk; \R^2)$ function defined as 
\begin{equation}
u(x,y) := \begin{cases}
\alp & {\rm if}~ (x,y) \in \Aa,
\\
\bet & {\rm if}~ (x,y) \in \Bb,
\\
\gam & {\rm if}~ (x,y) \in \Cc.
\end{cases}\label{def:u}
\end{equation}
Then
\begin{equation}\label{eq:main}
\A(u,\sourcedisk)\leq |\sourcedisk| + 
\min \Big\{
{\mathfrak A}_{12}(\Gamma) + 
{\mathfrak A}_{23}(\Gamma) + 
{\mathfrak A}_{31}(\Gamma)  : 
\elementXclip \in \Xc
\Big\}.
\end{equation}
\end{Theorem} 
This theorem says that the singular part of $\A(u,D)$ can be 
estimated from above by 
\begin{equation}\label{eq:infimum_G}
\inf \Big\{
{\mathfrak A}_{12}(\Gamma) + 
{\mathfrak A}_{23}(\Gamma) + 
{\mathfrak A}_{31}(\Gamma)  : 
\elementXclip \in \Xc
\Big\}
\end{equation}
and that such an infimum is a minimum. Intuitively, to 
``fill the hole'' in the graph of $u$ with smooth two-codimensional approximating
surfaces of graph-type, we start to imagine three minimal surfaces,
placed vertically over the jump of $u$, coupled together via a 
triple point in the target triangle $\targettriangle$ (notice that the union of these
three minimal surfaces, viewed in $D \times \R^2$, is not smooth in correspondence
of the source triple junction). Such a construction
depends on the choice of a target triple point, and on 
a connection $\Gamma$ passing through it, dictating
the boundary condition for the three minimal surfaces, over 
the sides of the triangle $\targettriangle$. Theorem \ref{teo:main}
asserts that the interesting part of the 
relaxed area of $u$, namely its singular part,  cannot be 
larger
than what we obtain by minimizing over all possible target triple points
and all corresponding 
connections. 
As a direct consequence of 
the results in \cite{BePa:10}, \cite{Sc:19},
when $u = u_{\rm symm}$ (and $\sourceorigin$ is the center of $D$),
the inequality in \eqref{eq:main} is an equality, and
the infimum in \eqref{eq:infimum_G} is achieved by the 
Steiner graph connecting the three vertices of $\targettriangle$ (the optimal triple point being the Steiner point, i.e., the barycenter of $\targettriangle$). 
This seems
to be an interesting result that could be stated purely as 
a problem of three entangled area-minimizing surfaces
(each of which lies in a half-space of $\R^4$, the three
half-spaces having only $\{0\}\times \R^2$ in common)
without referring to the relaxation of the functional 
$\area(\cdot,\sourcedisk)$.
We do not know whether, in general, the
Steiner graph  is still the solution of the minimization problem in \eqref{eq:infimum_G},
when no symmetry assumptions (the case we are considering here) 
 are required. 
However it is reasonable to expect that, if in the source 
we have symmetry, {\it i.e.}, the source triple junction is positioned at the center
of $D$ and $u$ 
jumps along three segments  meeting 
at equal $120^\circ$ angles, and if the target triangle $\targettriangle$ is close
to be equilateral, the inequality in \eqref{eq:main}
to be still an equality. In this respect, it is worthwhile to 
observe that showing a lower estimate,
for instance showing that, in certain cases, the inequality
in \eqref{eq:main} is an equality, seems difficult. 
One of the main technical 
obstructions is due to 
the poor control on the tangential derivative  of $v^\eps$
in proximity of the jump of a discontinuous $L^1$-limit function $v$
(see \cite{BePaTe:16}), where $\{v^\eps\}$
is a sequence of Lipschitz maps converging
in $L^1(\Omega; \R^2)$ to $v$, and satisfying the uniform bound
$\sup_\eps \area(v^\eps, \Omega) < +\infty$. 
We also notice that
the symmetrization 
methods of \cite{Sc:19} cannot be applied anymore, in view
of the lackness of symmetry. 

It is worth mentioning that the restriction that we assume
on the connections $\Gamma$, namely that each $\Gamma_i$ is a graph
(possibly with vertical parts) on the corresponding two sides of $\targettriangle$,
cannot be avoided in our approach: indeed, only under this graphicality
assumption we can solve the minimum problem in \eqref{eq:main} in the class
of surfaces which are graphs over the rectangles $\rettangolo_{ij}$. In turn,
the graphicality of such minimal surfaces allows to construct
the sequence $\{u^\eps\}$, see \eqref{def:uhstrip}. Removing the graphicality
assumption on $\Gamma$ requires some change of perspective, and needs
further investigation. 

The content of the paper is the following. 
In Section \ref{subsec:functions_of_bounded_variation} we recall
some properties of functions of bounded variation of one variable,
the definition of generalized graph (formula \eqref{eq:generalized_graph}), and the chain rule. 
In Section \ref{subsec:cartesian_currents} we recall some
properties of Cartesian currents carried by a BV-function.
 The functional $\AG$, appearing on the right hand side of \eqref{eq:main}, is introduced in Definition \ref{def:functional_G}. 
In Section \ref{sec:upperbound} we show that 
\begin{equation}\label{eqn:elhagrahwan}
\A(u, \sourcedisk) \leq |\sourcedisk| +\inf\big \{\AG(\elementXclip): \elementXclip\in X_\Lip\big \},
\end{equation} 
see Theorem \ref{thm:upper_bound_with_infimum}.
The proof is rather involved, mainly due to technical difficulties:
we first start by supposing that the jump of $u$ is piecewise linear
(Proposition \ref{prop:piecewise_linear}). Some work is required
to define $u^\eps$ on an $\eps$-strip around the jump
of $u$ and avoiding a neighbourhood of
the source triple junction 
(formula \eqref{def:uhstrip}) and to define $u_\eps$ in 
the missing neighbourhood of the source triple junction (step 3 
of the proof
of Proposition \ref{prop:piecewise_linear}): 
the construction must be done in such a way that $u^\eps$ remains
Lipschitz, and turns out to be rather involved in the three triangles $T_1^\eps, T_2^\eps, T_3^\eps$, see Figure \ref{fig:T_eps}.
In Section \ref{sec:AnalysisAG} we prove that the infimum in \eqref{eq:infimum_G}
is a minimum. The proof is achieved by defining 
a topology in the space $X$ which allows 
to prove the 
density of $X_\Lip$ in $X$ (Lemma \ref{Lem:density}), 
the continuity of the functional $\AG$ (Proposition \ref{Lem:AGcontinuous})
and the sequential compactness of $X$ (Theorem \ref{teo:closureXclip}). 
This latter result is also based on a uniform bound on the length
of the connections (Proposition \ref{pro:uniform_estimate_on_the_length}), which is 
a consequence of the graphicality assumptions on the connections. 
 
\section{Some preliminaries}\label{sec:notation}
In this section we 
recall some results on functions of bounded variation of one variable \cite{AmFuPa:00},
and on cartesian currents\cite{Giq}, needed in the sequel.

\subsection{Functions of bounded variation in the interval}\label{subsec:functions_of_bounded_variation}
Let $\I\subset \R$ be a bounded open interval and $\dirdatum \in \BV(\I)$; then  
\begin{itemize}
\item $\dirdatum$ is bounded, and it is
continuous up to an at most countable set of 
points of $\I$ denoted by $J_\dirdatum$ (jump set);
\item the right and left limits $\dirdatum(s_\pm)$ of $\dirdatum$ 
exist at any $s\in \I$; 
the right limit $\dirdatum(a_{+})$ 
and the left limit $\dirdatum(b_{-})$ exist.
Thus we may define $$\dirdatum_{+}(s):=\max\{\dirdatum(s_+), \dirdatum(s_-)\},\qquad \dirdatum_{-}(s):=\min\{\dirdatum(s_+), \dirdatum(s_-)\},\qquad s\in \I;$$
\item  the distributional derivative $\dirdatum^\prime$ of $\dirdatum$ splits as
$$\dirdatum^\prime=\dot\dirdatum ds+\dot \dirdatum^{ (j)}+\dot \dirdatum^{ (c)},$$
where $\dot\dirdatum ds$ is the absolutely continuous part and $\dot\dirdatum$ is the differential of $\dirdatum$ \cite[p.138 and Cor.3.33]{AmFuPa:00}, ${\dot\dirdatum}^{ (j)}$ and ${\dot\dirdatum}^{ (c)}$ are the jump and the Cantor part respectively. 
\end{itemize}
We shall always assume that $\dirdatum$ is a 
good representative in its $\Lone$ class such that 
$\dirdatum(s)=\dirdatum_{+}(s)$ for all $s\in \I$;
the pointwise variation of $\dirdatum$ 
is equal to the total variation $\vert \dirdatum ^\prime \vert(\I)$. 

The generalized graph of $\dirdatum$ is defined as
\begin{equation}\label{eq:generalized_graph}
\Gamma_\dirdatum:=\{(s,\theta \dirdatum(s_-)+(1-\theta)\dirdatum(s_+)): s\in \I,~ \theta \in [0,1]\}, 
\end{equation}
and the subgraph of $\dirdatum$ as 
\begin{equation*}
S\G_{\dirdatum,\I}:=\{(s,t)\in \I\times \R~ : ~t\leq \dirdatum(s)\}.
\end{equation*}
We recall that, if $\dirdatum\in\Lone(\I),$ then $\dirdatum\in \BV(\I)$ if and only if $S\G_{\dirdatum,\I}$ has finite perimeter in $\I \times \R$.
We denote by $\partial^-S\G_{\dirdatum,\I}$ the reduced boundary of $S\G_{\dirdatum,\I}$.

{We conventionally set $\dirdatum(a_-)= 0$, $\dirdatum(b_+)=0$; in this case we can define $\Gamma_\dirdatum$ as in \eqref{eq:generalized_graph} with $\I$ replaced by $\cI$, hence the generalized graph will always pass through the end points of the interval (with possibly vertical 
parts over $a$ and $b$).}

\smallskip

The following result can be found for instance in  \cite[p.486]{Giq}.
\begin{Theorem}\label{teo:chain_rule}
Let $g\in \C^1(\R)$ and $\dirdatum\in \BV(\I).$ Then 
$g\circ \dirdatum\in \BV(\I)$ and 
\begin{align*}
(g\circ \dirdatum)^\prime&
= g^\prime (\dirdatum) \dot\dirdatum ds+ g^\prime(\dirdatum) \dot\dirdatum^{(c)} 
\qquad {~\rm in~} 
\I \setminus J_\dirdatum 
\\
(g\circ\dirdatum)^\prime&
=\sum_{s\in J_\dirdatum}n(s,J_\dirdatum)
\Big[g(\dirdatum_+(s))-g(\dirdatum_-(s))\Big] 
\delta_s \qquad {~\rm in~}  
J_\dirdatum,
\end{align*}
where $n(s,J_\dirdatum)
:=\frac{\dirdatum(s_{+})-\dirdatum(s_{-})}{\vert \dirdatum(s_{+})-\dirdatum(s_{-})\vert}$ and 
$\delta_s$ is the Dirac delta at $s$.
\end{Theorem}

\subsection{Cartesian currents}\label{subsec:cartesian_currents}
Let $I\subset \R$ be a bounded open interval and $\dirdatum \in \BV(I)$. 
We denote by 
$$
[\![ S\G_{\dirdatum,I} ]\!]\in\mathcal D^2(\R^2)
$$
the 2-current in $\R \times \R$ defined as the integration over the subgraph $S\G_{\dirdatum,I}$. 
The current 
$[\![ S\G_{\dirdatum,I} ]\!]\res I\times \R$ 
can be also identified with an integer multiplicity current in $I\times\R$; 
moreover $S\G_{\dirdatum,I}$ has finite perimeter in $I\times \R$ so, if $\partial [\![ S\G_{\dirdatum,I} ]\!]\res I\times \R$  denotes the 1-current in $I\times \R$ defined as its boundary, this results of finite mass.

For future purposes we recall the following result, see \cite[Section 4.2.4]{Giq}. 

\begin{Theorem}\label{thm: boundary currents}
Let $\dirdatum \in \BV(I)$ and $\T$ be the current defined by
\begin{equation}
\T:= - \partial [\![ S\G_{\dirdatum,I} ]\!]\res I\times \R.
\end{equation}
Then $\T\in \mathcal D^1(I\times\R)$ 
 is a Cartesian current, and
\begin{equation}
 \T(\omega)= - \int <\omega(x), * \nu(x,S\G_{\dirdatum,I})> d\Hone \res \partial^- S\G_{\dirdatum,I}(x) \qquad \forall\omega\in \DD^1(I\times\R),
\end{equation}
where $*$ is the Hodge operator and $\nu(\cdot,S\G_{\dirdatum,I})$ is the inward generalized unit normal. Moreover $\T$ can be decomposed 
into three mutually singular currents 
\begin{equation}
\T=\T^{(a)}+ \T^{(j)}+\T^{(c)},
\end{equation}
such that
\begin{align}
\T^{(a)}(\omega)&=\int_I [ \omega_1(s,\dirdatum(s))+\omega_2(s,\dirdatum(s)) \dot\dirdatum (s) ]ds,\label{eqn:Ta}\\
\T^{(j)}(\omega)&=\sum_{s\in J_\dirdatum} n(s,J_\dirdatum) 
\int_{\dirdatum_{-}(s)}^{\dirdatum_{+}(s)} \omega_2(s,\sigma)d\sigma,
\label{eqn:Tj}\\
\T^{(c)}(\omega)&=\int_I \omega_2(s,\dirdatum(s)) \dot\dirdatum^{(c)},\label{eqn:Tc}
\end{align}
where $\omega=\omega_1ds+\omega_2d\sigma$.
\end{Theorem}

The current 
$\T$
is boundaryless in $I\times\R$, namely $\partial\T=0$. 
Furthermore, 
if $\Gamma_\varphi$ is the generalized graph of $\varphi$ 
as defined in \eqref{eq:generalized_graph}, it turns 
out that $\partial^-S\G_{\dirdatum,I}
\cap (I\times\R)$ is a subset of $\Gamma_\varphi$ and they 
differ of a $\mathcal H^1$-negligible set. Namely
$$\partial^-S\G_{\dirdatum,I}
\cap (I\times\R)\subseteq\Gamma_\varphi,\qquad \mathcal 
H^1(\Gamma_\varphi\setminus\partial^-S\G_{\dirdatum,I})=0. $$
It easily follows that the current $\T$ coincides with the integration over the rectifiable set $\Gamma_\varphi$ (with the correct orientation).

From now on, when the interval is clear from the context, we will simply denote $S\G_{\dirdatum,I}$ by $S\G_{\dirdatum}$.

%
\section{The functional $\mathcal G$}\label{sec:Connection}
In order to prove our main result (Theorem \ref{teo:main}) we need 
some preparation.
Take 
three open non-overlapping non-empty connected regions $ \Aa,~\Bb,~\Cc$ 
of an open disk $\sourcedisk$, each $E_i$ with non empty interior and with $\overline{\Aa} \cup \overline{\Bb} \cup \overline{\Cc}=\overline{\sourcedisk}$, and let {$\cij$} be their boundaries in $\sourcedisk$
as in the introduction.

Let $\alp,~ \bet, ~\gam$ be 
the vertices of a closed 
triangle $\targettriangle$ as in Section \ref{sec:introduction};
we suppose that 
$\targettriangle$ contains the origin $\targetorigin$ in its interior,
and let $\ell_{ij}$ be as in \eqref{eq:l_ij}.

\begin{Definition}[\textbf{Connections in $\targettriangle$}]\label{def:connections}

 We say that $\elementXclip:=(\Gammaalp, \Gammabet, \Gammagam)$ is a $\BV$ graph-type (resp. $\Lip$ graph-type) connection in $\targettriangle$ if  $\pGammaalpi, ~i\in \{1,2,3\},$ are subsets of $\targettriangle$ such that $\Gammaalp\cap \Gammabet=\Gammabet \cap \Gammagam=\Gammagam \cap \Gammaalp$ is one point $p$ of $\targettriangle$ called target triple point of $\Gamma$, $\alpi \in \pGammaalpi$ for any $i=1,2,3$, and  
$$\pGammaalpij:=\pGammaalpi \cup \pGammaalpj, \qquad ij\in \{12,23,31\},$$
 can be written as the generalized graph (resp. graph) of a function of bounded variation (resp. Lipschitz function) over the closed segment $\overline{\alpi \alpj}$  
(see Figure \ref{fig:targettriangle}).
 \end{Definition}
 
Note that the case $\genericp \in \partial \targettriangle$ is not excluded. However, by definition, 
if {$\pi_{{ij}}: \targettriangle \to \R_{\overline{\alpi \alpj}}$,} ${ij} \in \{12,23,31\}$, is the orthogonal projection on the line $\R_{\overline{\alpi \alpj}}$ containing $\overline{\alpi \alpj}$, then $\pi_{{ij}}(p)\in \overline{\alpi \alpj}$. Set 
 \begin{equation}\label{eqn:wij}
 \wij:=|\alpi-\pi_{{ij}}(\genericp)|.
 \end{equation}

If necessary, in the sequel we will often identify $\pGammaalpij$ with the (generalized) graph $\Gamma_{\phiij}$ of a function 
\begin{align}\label{varphi_ij}
 \phialpij : [0,\ellalpij] \to [0, {\rm diam} \targettriangle], \qquad \phialpij=\phialpij(\pGammaalpij),
\end{align}
 of bounded variation. If $\targettriangle$ is acute, choosing a suitable cartesian coordinate system where the $s$-axis is the line  $\R_{\overline{\alpi \alpj}}$, we necessarily have  $\phiij(0)=\phiij(\ellalpij)=0$. In contrast, if the angle of $\targettriangle$ at $\alpha_i$ is greater than or equal to $\frac{\pi}{2}$ then $\varphi_{ij}$ might have a vertical part over $\alpha_i$ and  $\phiij(0_+)>0$. In this case the generalized graph of $\varphi_{ij}$ does not pass through $\alpha_i$.
 
 In the sequel it will be often convenient to consider an extension  of $\varphi_{ij}$ on $(-\infty,0)\cup(\ellalpij,+\infty)$. This extension is denoted by $\tilde\varphi_{ij}$. In the case of acute triangle $\tilde \varphi_{ij}$ is always set equal to $0$ on $(-\infty,0)\cup(\ellalpij,+\infty)$.
 
 %
 \begin{Remark}\label{Rmk:graph.} \rm
 If for any ${ij} \in \{12,23,31\}$, $\wij$ in \eqref{eqn:wij} is a point of 
continuity of $\phialpij$ then the intersection of the generalized graph of $\phialpki$ and the set $[\wki,\ellalpki]\times \R$ coincides with $\Gamma_i$ which is also the intersection of the generalized graph of $\phialpij$ 
with the set $[0,\wij]\times \R$, where ${ij} ,\ki \in \{12,23,31\},~{ij} \not =\ki$. If $\wij$ is a 
discontinuity point of $\phiij$ this is in general not true, as in Figure \ref{fig:caseii}, when $i=2$.
 
 \end{Remark} 
 \begin{Remark}\label{H} \rm
 {Assume that an angle of $\targettriangle$ is greater than $\frac{\pi}{2}$, say for instance the angle at $\alpha_1$; as already said,  the generalized graphs composing a connection $\Gamma$   are allowed to have vertical parts over $\alpha_1$. The target triple point $p$ of any connection $\Gamma$ belongs to $T_{int}\subset T$, the part of the triangle $\targettriangle$ which is enclosed between the two lines passing through $\alp$ and orthogonal to $\overline{\alp\bet}$ and $\overline{\alp \gam}$ respectively.}
 \end{Remark}

  Define the classes:
  \begin{align}
 X_{\Lip}&:=  \big\{ \elementXclip:  \elementXclip{~} \Lip {\rm ~graph -type ~connection ~ in~} \targettriangle 
   \big\},\label{eqn:classXlip}\\
   X&:=  \big\{\elementXclip:  \elementXclip{~} \BV {\rm ~graph -type ~connection~ in~} \targettriangle \big\}.\label{eqn:classX}
\end{align}
Obviously $X_{\Lip}\subset X$.

 \subsection{Useful results on one-codimensional area-minimizing cartesian surfaces}\label{tenda}
 Let
$\rettangolo_{ij}$ be as in \eqref{eq:R_ij}, and $\elementXclip \in X$. Then $\pGammaalpij,~ij\in\{12,23,31\} $ are (generalized) graphs of functions $\phiij$ of bounded variation over  $[0,\ellalpij]$. Let $\Ball \subset \R^2$ be  an open disk containing the doubled rectangle $\widehat\rettangolo_{i j}$ 
defined as
\begin{equation}
\widehat\rettangolo_{i j}:=\left[0, \ellalpij \right] \times [0,2\radiusalpij],\qquad ij\in\{12,23,31\} .\label{eqn:Rhat}
\end{equation}
We use for simplicity the same notation $\phialpij$ for the 
extension of $\phialpij$ to $\widehat \rectangle_{i j}$,
 defined as 
 \begin{equation}\label{eq:extension_of_varphi_bis}
 \phialpij (s,t)=\phialpij(s), 
\qquad (s,t)\in \widehat\rettangolo_{i j}  
\qquad ij\in\{12,23,31\}\},
 \end{equation}
and for the extension of $\phialpij$ to a $W^{1,1}$ function on $\Ball \setminus \widehat\rettangolo_{i j}$ as in \cite[Theorem 2.16]{Gi:84}. 

Let $\widehat\tenda_{i j}=\widehat\tenda_{i j}(\Gamma)$,  $ij\in\{12,23,31\}$, be the solution of following Dirichlet minimum problem:
\begin{align}\label{widehat_tenda_pen}
\min 
\left\{
\int_{\widehat\rettangolo_{i j}}
\sqrt{1 + \vert D f\vert^2} ~+\int_ {\partial \widehat \rettangolo_{i j}}\vert f-\phialpij\vert d \Hone: ~ f \in \BV (\Ball), 
~ f = \phialpij ~{\rm on}~ 
\Ball \setminus \widehat\rettangolo_{i j}\right\},
\end{align}
where {$\int_{\widehat\rettangolo_{i j}}
\sqrt{1 + \vert D f\vert^2}$ is the extension of the area functional to $\BV(\widehat\rettangolo_{i j})$ as defined in \cite[Definition 14.1]{Gi:84}. }

From \cite[Theorem 15.9]{Gi:84}  and the fact that the restriction of $\phialpij $ to ${\partial  \widehat\rettangolo_{i j}}$ is continuous up to a countable set of points, it follows that $\widehat \tenda_{i j}$ solves also
\begin{align}\label{widehat_tenda}
\min 
\left\{
\int_{\widehat \rettangolo_{i j}}
\sqrt{1 + \vert \grad f\vert^2} ~ds dt : ~ f \in W^{1,1}(\widehat\rettangolo_{i j}), 
~ f = \phialpij ~\H^1 -a.e.~{\rm on}~ 
\partial \widehat \rettangolo_{i j}
\right\}.
\end{align}
Let $\tenda_{i j}=\tenda_{i j}(\Gamma)$ be the restriction of $\widehat \tenda_{i j}$ to $\rettangolo_{i j}$. Then, by the symmetry of $\phialpij$ with respect to the line $\{t=\radiusalpij\}$, $\tenda_{i j}$ is the unique solution of the Dirichlet-Neumann minimum problem
\eqref{eq:Dirichlet_Neumann_minimum_problem}.
From \eqref{eq:extension_of_varphi} it follows that the Dirichlet condition is zero on the sides $\{0\}\times [0,\radiusalpij]$ and $ \{\ellalpij\}  \times [0,\radiusalpij] $ of the rectangle $\rettangolo_{i j}$. Note that $\tenda_{i j}$ is analytic in the interior of $\rettangolo_{i j}$ {but {\it not} necessarily Lipschitz in $\rettangolo_{i j}$ \cite[Theorem 14.13]{Gi:84}, see Figure \ref{fig:graph_m_i}.}

\begin{Definition}[\textbf{The functional $\AG$}]\label{def:functional_G}
We define the functional
$\AG:X \longrightarrow [0,+\infty)$
as
\begin{equation}\label{eqn:AG}
\AG(\elementXclip):={\mathfrak A}_{12}(\Gamma)
 + {\mathfrak A}_{23}(\Gamma) +{\mathfrak A}_{31}(\Gamma),
\end{equation}
{where ${\mathfrak A}_{ij}(\Gamma)$ are as in \eqref{eq:A_ij}.}
 \end{Definition}
The properties of the functional $\AG$ will be discussed in Section \ref{sec:AnalysisAG}.

\section{Infimum of $\AG$ as an upper bound of $\A(u, D)$} \label{sec:upperbound}
The aim of this section is to provide
the following upper bound for 
$\A(u,\sourcedisk)$.

\begin{Theorem} \label{thm:upper_bound_with_infimum}
Let $u\in BV(\sourcedisk; \{\alp, \bet, \gam \})$ be the function 
defined in \eqref{def:u}. Then 

\begin{equation}\label{Gammalimsupcurves}
\A(u, \sourcedisk) \leq |\sourcedisk| +\inf\big \{\AG(\elementXclip): \elementXclip\in X_\Lip\big \}.
\end{equation} 
\end{Theorem}

It is not difficult to see, 
by truncating the minimal surfaces in \eqref{eq:Dirichlet_Neumann_minimum_problem} with the lateral boundary of the prisms $[0,\ellalpij] \times \targettriangle$, that the infimum in \eqref{Gammalimsupcurves} is the same as the infimum obtained without requiring in Definition \ref{def:connections} that $\pGammaalpi \subset \targettriangle$, $i\in \{1,2,3\}$.

 \begin{Lemma}\label{lem:phiapprox}
Let $\ell\geq 0,~\mathrm{\genericp}\geq 0,~\dirdatum \in \Lip([0,\ell];[0,+\infty))$ be such that $\dirdatum(0)=\dirdatum(\ell)=0$ and $\wo,\in [0,\ell]$ so that $\dirdatum(\wo)={\mathrm{ \genericp}}$. Then there exists a sequence $\{\dirdatum^{\sigma}\}$ of $\C^{\infty}$ equi-Lipschitz functions in $[0,\ell]$, converging to $\dirdatum$ in $\Lone ([0,\ell])$ and uniformly on $[0,\ell]$ as $\sigma \to 0^+$, such that $$\dirdatum^{\sigma}(0)=\dirdatum^{\sigma}(\ell)=0,\qquad \dirdatum^{\sigma}(\wo)=\mathrm{ \genericp}, \qquad{\rm for~ any~ }\sigma>0.$$
\end{Lemma}

\begin{proof}
Let us extend $\dirdatum$ in $\R$ such that $\dirdatum(s)=0$ in $\R \setminus [0,\ell]$, so that the extension (still denoted by $\dirdatum$) belongs to $\Lip(\R)$. Let 
 $\widehatphialpbet^{\sigma}(s):=\eta_\sigma * \dirdatum$ in $ \R$,
 where $\{\eta_{\sigma}\}$ is a standard sequence of mollifiers. Hence $\widehatphialpbet^{\sigma}\in \C^\infty(\R)$, $\Lip (\widehatphialpbet^{\sigma})\leq \Lip(\dirdatum)$ and the sequence $\{\widehatphialpbet^{\sigma}\}$ converges uniformly to $\dirdatum$ on compact subsets of $\R$. Without 
loss of generality we may assume 
$\widehatphialpbet^{\sigma}(s)=0$ in $\R\setminus (-\sigma/2,\ell+\sigma/2)$ and $\widehatphialpbet^{\sigma}(\frac{\ell+2\sigma}{\ell}\wo-\sigma)=\mathrm{ \genericp}+c_\sigma,~ c_\sigma=o(1)$. {Let us first suppose $\mathrm{ \genericp}\neq0$}. We  define 
 \begin{gather}
  \dirdatum^{\sigma}: [0,\ell]\to [0,+\infty), \qquad
\dirdatum^{\sigma}(s):=\frac{\mathrm{ \genericp}}{\mathrm{ \genericp}+c_\sigma}~\widehatphialpbet^{\sigma}\big(\frac{\ell+2\sigma}{\ell}~s-\sigma\big).\label{eqn:phisigma}
  \end{gather}
  It is easy to see that
$\dirdatum^{\sigma} \in \C^{\infty}([0,\ell])$, 
$\dirdatum^{\sigma} (0)= \dirdatum^{\sigma}(\ell)=0,~ \dirdatum^{\sigma}(\wo)=\mathrm{ \genericp},$ 
$\dirdatum^{\sigma}$ are equi-Lipschitz, and $\{\dirdatum^{\sigma}\}$ converges to $\dirdatum$ in $L^1([0,\ell])$ as $\sigma \to 0$. Notice that the obtained approximation is constantly null in a neighborhood of $0$ and $\ell$. 

In the case $\mathrm{ \genericp}=0$, we argue differently. We consider the two intervals $[0,\wo]$ and $[\wo,\ell]$ and we repeat the same approximation above in the single intervals; more precisely we choose two points $w_1\in (0,\wo)$ and $w_2\in(\wo,\ell)$ with $\varphi(w_1)>0$, $\varphi(w_2)>0$ (if these points does not exist it means that the functions are constantly $0$ and they are already smooth, so there is nothing to prove). Then we approximate the two functions $\varphi\res (0,\wo)$ and $\varphi\res [\wo,\ell]$ as before, and we glue them along $w$. Note that the glued function is smooth in $w$ since both the two smooth approximations are constantly $0$ in a neighborhood of $\wo$. 

\end{proof}

To prove Theorem \ref{thm:upper_bound_with_infimum} we 
use the three area-minimizing functions $\tenda_{ij}, {ij}\in\{12,23,31\},$ introduced in Section \ref{tenda}, to construct a sequence $\{\uh\}$ of Lipschitz functions that converge to $u$ in $\Lone(\sourcedisk;\R^2)$. However $\tenda_{ij}, {ij}\in\{12,23,31\},$ are only locally Lipschitz so we need the following smoothing lemma.

\begin{Lemma}\label{lem:msigma}
Let $\elementXclip \in X_\Lip$, ${ij} \in \{12,23,31\}$. Let $\phiij =\phiij(\Gamma_{ij})\in \Lip([0,\ellalpij])$, $\tenda_{i j}=\tenda_{i j}(\Gamma_{ij})\in W^{1,1}(\rettangolo_{ij}),$ be defined as in Section \ref{tenda}. Then there exists a sequence $\{\tenda_{ij}^{\sigma}\}$ of Lipschitz functions such that $\tenda_{ij}^{\sigma}:\rettangolo_{ij}\to \R$, $\tenda_{{ij}}^{\sigma}=\varphi_{ij} {~\rm on~}{\partial_D \rettangolo_{ij}},$ and 
 \begin{equation}\label{lip approx m}
\left\vert 
\int_{\rettangolo_{ij}} \sqrt{1+ 
\vert \grad \tenda_{ij}\vert^2}~dsdt 
- 
\int_{\rettangolo_{ij}}\sqrt{1+ \vert \grad \tenda_{ij}^{\sigma}\vert^2}
~dsdt\right\vert \leq O(\sigma).
\end{equation}
\end{Lemma}
\begin {proof}
This can be easily proved using an argument similar to the one in \cite[p.378: p.381]{BePa:10}, and using also Lemma \ref{lem:phiapprox} with the choice $\ww=\wij$ and $\mathrm{ \genericp}=\varphi_{ij}(w_{ij})$. \end{proof}

We start to prove Theorem \ref{thm:upper_bound_with_infimum} in the special case of a piecewise linear jump, as in Figure \ref{fig:domain_segment_1}.

\begin{Proposition}\label{prop:piecewise_linear}
Let $u\in BV(\sourcedisk; \{\alp, \bet, \gam \})$ be the map 
defined in \eqref{def:u} and assume that the jump set of $u$ consists of three distinct segments that meet at the origin and reach the boundary of $\sourcedisk$. Then \eqref{Gammalimsupcurves} holds.
\end{Proposition}
\begin{proof}
Let $\elementXclip \in X_\Lip$ be a connection passing through $\genericp \in \targettriangle$  and $\AG(\elementXclip):={\mathfrak A}_{12} (\elementXclip)+ {\mathfrak A}_{23}(\elementXclip) +{\mathfrak A}_{31}(\elementXclip). $
To prove the proposition it is sufficient to construct a sequence $\{ \uh \} \subset \Lip(\sourcedisk; \R^2)$ converging to $u$ in $L^1(\sourcedisk;\R^2)$ such that
 \begin{equation}
 \lim_{\eps\to 0}\A(\uh,\sourcedisk) \leq |\sourcedisk| +{\mathfrak A}_{12}(\elementXclip)+{\mathfrak A}_{23}(\elementXclip)+{\mathfrak A}_{31}(\elementXclip).\label{eqn:case1f}
 \end{equation}

{\it Case 1.} Assume that the segments separating $\Aa,~ \Bb,~ \Cc$ meet at the origin with angles less than $\pi$, as in Figure \ref{fig:domain_segment_1}. 

\begin{figure}
\centering
  \begin{center}
    \begin{tikzpicture}
        \draw[black,thin](0,0)--(3.49,0)node[black,below=8,left=35]{$\radiusbetgam$};
         \draw[white,thick](1.1,1.1)circle(1pt) node[black,right=2]{$\Bb$};
          \draw[black,thin](0,0)--(0,2.49)node[black,left=8,below=35]{$\radiusalpbet$};
         \draw[white,thick](-1.8,0.5)circle(1pt) node[black,right=2]{$\Aa$};
         \draw[black,thin](0,0)--(-2.07,-2.07)node[black,right=25,above=27]{$\radiusgamalp$};
      \draw[white,thick](0.5,-1.5)circle(1pt) node[black,right=2]{$\Cc$};
         \draw[black,thick](0.5,-0.5)circle(3);   
        \filldraw[black!80](0,0)circle(1pt) node[black,below=7,right]{$0_S$};
    \end{tikzpicture}      
  \end{center}
 \caption{$\Aa,~\Bb,~\Cc$ are separated by three segments of length $\radiusalpbet$, $\radiusbetgam$, $\radiusgamalp$ that meet at the origin.}
 \label{fig:domain_segment_1}
\end{figure}
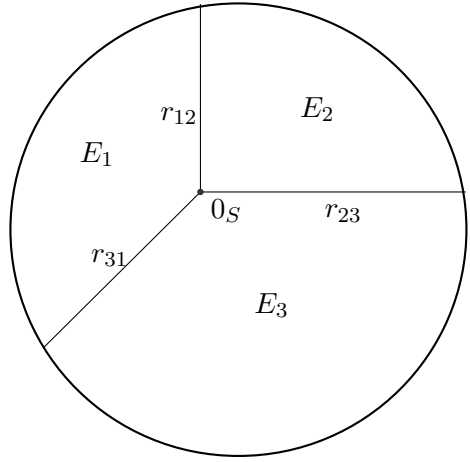

 To simplify the computation we may assume that
$p=\targetorigin,$ see Figure \ref{fig:targettriangle}. The idea 
of the proof is similar to the one used in \cite{BePa:10}, with however new difficulties, in particular in $T^\eps$ (step 3). We will specify various subsets of $\sourcedisk$ and define the sequence $\{\uh\}$ on each of these sets. 
Let $\eps >0$ be sufficient small and $\del_\eps >0$ be such that $\del_\eps \to 0$ as $\eps \to 0$. Define $\Triangle^\eps$ to be the triangle with the  origin $\sourceorigin$ in its interior, with {vertices $\zeta^1=\zeta^1_\eps$, $\zeta^2=\zeta^2_\eps$, and $\zeta^3=\zeta^3_\eps$, and} sides of lengths $\eps_\cc, ~\eps_\ac,~\eps_\bc$, {$\eps_{ij}:=|\zeta^i-\zeta^j|$}; {the sides of $\Triangle^\eps$ are perpendicular to the lines containing $r_{12}$, $r_{23}$, $r_{31}$ (respectively) and their distance from the origin $\sourceorigin$ equals $\del_\eps$}. Define three cygar-shaped sets $ S_\ac ^\eps,~S_\bc ^\eps$ and $S_\cc ^\eps$ as in Figure \ref{fig:Subsets_Br}: if for instance $y$ is a coordinate on $r_{12}$ and $x$ is the perpendicular coordinate, then $S_\cc^\eps$ is defined as
\begin{equation}\label{def:cygar}
S_{\cc}^\eps 
:= \left\{(x,y) \in \sourcedisk : x\in (\etia, \etiia),~
 y \geq \del_\eps  
 \right\},
\end{equation}
where $$\zeta^i=(\zeta^i_1,\zeta^i_2), \;\;\;i=1,2,3.$$ 

Let us set
\begin{equation}
\Aae
:=
\Aa \setminus \left(S_\bc^\eps \cup T^\eps \cup S_\cc^\eps \right),
\quad
\Bbe := 
\Bb \setminus \left(S_\ac^\eps \cup T^\eps \cup S_\cc^\eps \right), 
\quad \Cce := 
\Cc \setminus \left(S_\ac^\eps \cup T^\eps \cup S_\bc^\eps \right).\label{eqn:4.51}
\end{equation}
\smallskip
{\it Step 1.} Definition of $\uh$ on $\Aae \cup \Bbe \cup \Cce$.
We define
\begin{equation}\label{def:uhregioni}
\uh := 
\begin{cases}
\alp & {\rm in}~ \Aae,
\\
\bet & {\rm in}~ \Bbe,
\\
\gam & {\rm in}~ \Cce.
\end{cases}
\end{equation}
Note that 
$\A(\uh, \Aae \cup \Bbe \cup \Cce)  = 
\vert \Aae
\vert +
\vert \Bbe 
\vert + \vert \Cce 
\vert$, hence
\begin{equation}\label{areanulla}
\lim_{\eps \to 0^+}
\A(\uh, \Aae \cup \Bbe \cup \Cce) = \vert \sourcedisk\vert.
\end{equation}
\begin{figure}
\centering
\begin{subfigure}{.4\textwidth}
  \centering
  \begin{center}
    \begin{tikzpicture}
        \draw[black!50,thin,dashed](0,0)--(4.46,0)node[black,below=13,left=2]{$S_\ac^\eps$};
        \draw[black](0.4,4)--(0.4,-1)--(4.40,-1);
        
                 \draw[white](1.7,2)circle(1pt) node[black,right=2]{$\Bbe$};
          \draw[black!50,thin,dashed](0,0)--(0,4.0)node[black,below=15,left=2]{$S_\cc^\eps$};
                  \draw[black](4.5,0.4)--(-1,0.4)--(-1,3.7);
                  
         \draw[white](-1.8,0.5)circle(1pt) node[black,right=4]{$\eti$};
\draw[black](0.4,0.4)circle(1pt) node[black,right=7,above]{$\etii$};
\draw[white](0.42,-1.02)circle(1pt) node[black,below=7,right]{$\etiii$};
          
         \draw[white](-3,1)circle(1pt) node[black,right=2]{$\Aae$};
         
         \draw[black!50,thin,dashed](0,0)--(-2.58,-2.58)node[black,above=10]{$S_\bc^\eps$};
         \draw[black](-3.08,-1.7)--(-1.0,0.4)--(-0.3,-0.3)--(0.42,-1.02)--(-1.8,-3.28);

      \draw[white](1.2,-2.3)circle(1pt) node[black,right=2]{$\Cce$};
         
         \draw[black,thick](0.5,0)circle(4);
         
        \filldraw[white!80](0.1,-0)circle(1pt) node[black,left=0.001]{$T^\eps$};

    \end{tikzpicture}      
  \end{center}
 \caption{Case 1 of the proof of Proposition \ref{prop:piecewise_linear}}
  \label{fig:Subsets_Br}
\end{subfigure}%
\begin{subfigure}{.6\textwidth}
  \centering
  \begin{center}
    \begin{tikzpicture}
       
       \draw[black](-2,1)--(2,1)--(2,-3)--(-2,1);
       \draw[black](-0.5,1)--(2,-0.5)--(0,-1)--(-0.5,1);
       \draw[black](-2,1)circle(1pt) node[black,left=4]{$\eti$};
\draw[black](2,1)circle(1pt) node[black,right=7,above]{$\etii$};
\draw[black](2,-3)circle(1pt) node[black,below=7,right]{$\etiii$};

       \draw[black](-0.5,1)circle(1pt) node[black,above]{$w^a$};
       \draw[black](2,-0.5)circle(1pt) node[black,right]{$w^b$};
       \draw[black](0,-1)circle(1pt) node[black,left=2]{$w^c$};
       
       \draw[white](-0.8,0.5)circle(1pt) node[black]{$T^\eps_1$};
       \draw[white](1,0.5)circle(1pt) node[black]{$T^\eps_2$};
       \draw[white](1.3,-1.5)circle(1pt) node[black]{$T^\eps_3$};
       \draw[white](0,0,4)circle(1pt) node[black]{$T^\eps$};
       \draw[white](0.4,-0.2)circle(1pt) node[black]{$T^\eps_0$};

    \end{tikzpicture}       
  \end{center}
  \caption{Zoom of $T^\eps$ in (a).}
  \label{fig:T_eps}
\end{subfigure}
\caption{}
\label{fig:division}
\end{figure}

\smallskip

{\it Step 2.} Definition of $\uh$ on $S^\eps_\ac \cup S^\eps_\bc \cup S^\eps_\cc$ .

We will start with the construction on $S^\eps_\cc$. Set 
$$
\xi = (\xi_1,\xi_2) := \frac{\bet - \alp}{\ell_\cc} \in \mathbb S^1, \qquad 
\eta = (\eta_1, \eta_2) 
:= \xi^\perp,
$$
where $^\perp$  denotes the counterclockwise rotation of $\pi/2$. 

Let $\psi_\cc^\eps : 
\left[\del_\eps,\raggio_\cc+\ce \right]
 \to
 [0,\raggio_\cc]$
be linear, increasing, surjective, where $\ce>0$ is the smallest number such that 
$$S^\eps_\cc \subset[\etia,\etiia]\times\left[\del_\eps,\raggio_\cc+\ce \right],\qquad  \lim_{\eps \to 0^+}\ce=0 .$$ 
Note that for any $y \in \left[\del_\eps,\raggio_\cc+\ce\right]$ we have
\begin{equation}\label{kappaeps}
(\psi_\cc^\eps)'(y) = 
\frac{\raggio_\cc}{\raggio_\cc+\ce
- \del} =: \kappa_\eps,  \qquad 
\lim_{\eps \to 0^+} \kappa_\eps=1. 
\end{equation}
Let $ \tenda_\cc^{\sigma}$ be the map defined in Lemma \ref{lem:msigma}, whose area on $\rettangolo_{12}$ is by construction close to ${\mathfrak A}_{12}$, with $\{\sigma_\eps\}\subset (0,+\infty)$ a sequence such that
\begin{equation}\label{chegabbarmivolete}
\lim_{\eps \to 0^+} \sigma_\eps =0.
\end{equation}

We set, with $\sigma=\sigma_\eps$ for simplicity, 
\begin{equation}\label{def:uhstrip}
\uhs(x,y) := \alp + \left( \frac{x-\etia}{\eps_\cc}
\right)\ell_\cc  \xi + \tenda_\cc^{\sigma}\left(\frac{ x-\etia}{\eps_\cc}\ell_\cc
~,~\psi_\cc^\eps(y)\right) \eta, \qquad (x,y) \in S^\eps_\cc.
\end{equation}
Observe that $\uhs =(\uhs_1,\uhs_2) \in \Lip(S^\eps_\cc; \R^2)$,
$\uhs =\alp$ on $\{(x,y) \in S^\eps_\cc : x = \etia\}$, 
and $\uhs = \bet$ on $\{(x,y) \in S^\eps_\cc : x = \etiia\}$. By the definition of $\tenda_\cc^{\sigma},$ it is uniquely defined the point (depending on $\eps$) $w^a=(w^a_1,w^a_2)\in \overline{\eti \etii}$ such that  $\uhs(w^a_1,w^a_2)=0_T$ (see Figure \ref{fig:T_eps}).
Write for simplicity 
$$
\widetilde \tenda = 
\tenda_\cc^{\sigma}.
$$

Using that $\vert \xi\vert= \vert \eta\vert =1$, $\xi_1 \eta_1 + \xi_2 \eta_2=0$, and $\xi_1\eta_2 - \xi_2 \eta_1= 1$, we compute
$$ 1 + \vert \grad \uhs_1\vert^2 + \vert
\grad \uhs_2\vert^2 + \left(\frac{\partial \uhs_1}{\partial x}
\frac{\partial \uhs_2}{\partial y}-\frac{\partial \uhs_1}{\partial y}
\frac{\partial \uhs_2}{\partial x}\right)^2
=  
1 + 
\frac{\ell_\cc^2}{\eps_\cc^2} \left( 
1+ \big(
\widetilde \tenda_s \big)^2 
+ 
\big(
\widetilde \tenda_t
\big)^2
\kappa_\eps^2 \left(1+
\frac{\eps_\cc^2}{\ell_\cc^2}
\right)
\right),
$$
where 
$\widetilde \tenda_s, \widetilde \tenda_t$ denote, respectively, the partial derivatives
of $\widetilde \tenda$ with respect to $s:=\frac{ x-\etia}{\eps_\cc}\ell_\cc$ and $t:=\psi_\cc^\eps(y)$, and 
are evaluated
at $\left(\frac{ x-\etia}{\eps_\cc}\ell_\cc~,~\psi_\cc^\eps(y)\right)$.
As a consequence
\begin{align}\label{contributonellastriscia}
& \A(\uh, S^\eps_\cc)\nonumber
\\ =& \frac{\ell_\cc}{\eps_\cc}
 \int_{S^\eps_\cc}
\sqrt{
1 + \left[\widetilde \tenda_s\left(\frac{ x-\etia}{\eps_\cc}\ell_\cc,~\psi_\cc^\eps(y)\right)\right]^2
+ 
\left[\widetilde \tenda_t\left(\frac{ x-\etia}{\eps_\cc}\ell_\cc,~\psi_\cc^\eps(y)\right)\right]^2
\kappa_\eps^2 
\left(1 + \frac{\eps_\cc^2}{\ell_\cc^2}
\right)
+ O(\eps^2)
}~dxdy\nonumber
\\
= &  
\frac{1}{\kappa_\eps} \int_{\rettangolo_\cc\setminus P_\eps} 
\sqrt{
1 + \left[\widetilde \tenda_s\left(s,t\right)\right]^2 
+ 
\left[\widetilde \tenda_t\left(s,t\right)\right]^2
\kappa_\eps^2\left(1 + \frac{\eps_\cc^2}{\ell_\cc^2}
\right)
 + O(\eps^2)
}
~dsdt,
\end{align}
where the last equality follows by the change of variables
$$
\Phi:  
\rettangolo_\cc\ni (s,t)\mapsto 
\Phi(s,t) := 
\left(\frac{\eps_\cc}{\ell_\cc}s +\etia, {\psi^\eps_\cc}^{-1}(t)\right) = (x,y)
\in \left[\etia, \etiia\right] \times 
\left [\del_\eps,\raggio_\cc+\ce\right] 
\supset S^\eps_\cc,
$$
and $P_\eps := \rettangolo_\cc \setminus \Phi^{-1}(S^\eps_\cc)$
(see Figure \ref{fig:Peps}).
\begin{figure}
\begin{center}
\includegraphics[height=3.0cm]{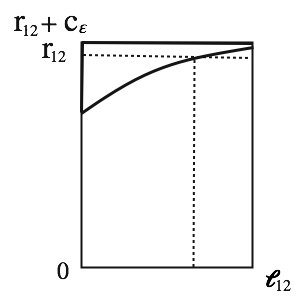}
\caption{\small{The set $P_\eps$ is bounded by
the bold contour.
}}
\label{fig:Peps}
\end{center}
\end{figure}
Hence, recalling also  \eqref{kappaeps}, we conclude 
\begin{equation}\label{contributonellastriscialimite}
\lim_{\eps \to 0^+} 
\A(\uhs, S^\eps_\cc) = 
 \int_{\rettangolo_\cc} 
\sqrt{
1 + 
\big(\widetilde \tenda_s
\big)^2
+ 
\big(\widetilde \tenda_t
\big)^2
}
~dsdt.
\end{equation}
We recall that from \eqref{lip approx m} it follows that
\begin{equation}\label{isitimmediate}
 \int_{\rettangolo_\cc} 
\sqrt{
1 + 
\big(\widetilde \tenda_s
\big)^2
+ 
\big(\widetilde \tenda_t
\big)^2
}
~dsdt = 
{\mathfrak A}_\cc(\Gamma) + O(\eps).
\end{equation}
Hence,  employing the same construction in the strips $S_\ac^\eps$  and $S^\eps_\bc$ we obtain
\begin{equation}\label{3cygar_contribution}
\lim_{\eps \to 0^+} 
\A(\uhs, 
S^\eps_\ac \cup 
S^\eps_\bc \cup 
S^\eps_\cc) =  {\mathfrak A}_\cc(\Gamma) +{\mathfrak A}_\ac(\Gamma) + {\mathfrak A}_\bc(\Gamma).
\end{equation}
\smallskip

{\it Step 3.} Definition of $\uh$ on $T^\eps$. 
We divide $T^\eps$ into four closed triangles $T^\eps_1$,
$T^\eps_2$, $T^\eps_3$ and $T^\eps_0$  as in  Figure \ref{fig:T_eps}. 
We set 
\begin{equation}\label{T0_contribution}
\uh := \targetorigin \qquad {\rm in}~ T^\eps_0.
\end{equation}
We first define $\uh$ on $\partial T^\eps_1$ as follows:
\begin{itemize}
\item[(i)] the value of $\uh$ at $\eti$ 
is $\alp$;
\item[(ii)] the value of $\uh$ on the side $\overline{w^cw^a}$ is 
 $\targetorigin$.
\end{itemize}
{Note that $\uh$ is already defined on the edges $\overline{\zeta^1w^a}$ and $\overline{\zeta^1w^c}$ and its graph over both edges is given by a rescaled version of the curve $\Gamma_1$ suitably parametrized.}

More precisely, we recall that 
 $\pi_{12}:\Gamma_1\rightarrow \overline{\alpha_1\alpha_2}$ and $\pi_{31}:\Gamma_1\rightarrow \overline{\alpha_3\alpha_1}$ are the orthogonal projections onto the edges $\overline{\alpha_1\alpha_2}$ and $\overline{\alpha_3\alpha_1}$. Since $\Gamma_1$ is, by hypothesis, a part of a Lipschitz graph, the maps $\pi_{12}\res\Gamma_1$ and $\pi_{31}\res\Gamma_1$ are bi-Lipschitz bijections between $\Gamma_1$ and the segments $\overline{\alpha_1\pi_{12}(p)}$ and $\overline{\alpha_1\pi_{31}(p)}$, respectively. We know that if $(s,t)$ are coordinates on $\targettriangle$ with respect to the system with  $s$-axis $\overline{\alpha_1\alpha_2}$, then the inverse of $\pi_{12}\res\Gamma_1$ is given by $\Phi_{12}:\overline{\alpha_1\alpha_2}\rightarrow\Gamma_1$,
$$\Phi_{12}((s,0))= (s,\varphi_{12}(s)).$$
Let us denote by $L_{12}=L^\eps_{12}:\overline{\zeta^1w^a}\subset \R_S^2\rightarrow \overline{\alpha_1\pi_{12}(p)}\subset\R_T^2$ and $L_{31}=L^\eps_{31}:\overline{\zeta^1w^a}\subset\R_S^2\rightarrow \overline{\alpha_1\pi_{31}(p)}\subset\R_T^2$ the linear bijective maps
$$L_{12}(Q)=\alpha_1+\frac{|Q-\zeta^1|}{\eps_{12}}(\alpha_2-\alpha_1),\qquad L_{31}(Q')=\alpha_1+\frac{|Q'-\zeta^1|}{\eps_{31}}(\alpha_3-\alpha_1).$$
Then we define
\begin{align}
 \uh:=\Phi_{12}\circ L_{12}\qquad \text{ on }\quad \overline{\zeta^1w^a},
\end{align}
and
\begin{align}
 \uh=\Phi_{31}\circ L_{31}\qquad \text{ on }\quad \overline{\zeta^1w^c}.
\end{align}
compare formula \eqref{def:uhstrip}. 
Since $\Phi_{12}$ and $\Phi_{31}$ are Lipschitz with Lipschitz constant independent of $\eps$, and the Lipschitz constants of $L_{12}$ and $L_{31}$ have order $\frac{1}{\eps}$, it follows that the Lipschitz constants of $\uh$ over the segments $\overline{\zeta^1w^a}$ and $\overline{\zeta^1w^c}$ have order $\frac{1}{\eps}$.

Now we want to define $\uh$ in the interior of $T^\eps_1$. First we observe that the map $\pi_{31}\circ \Phi_{12}:\overline{\alpha_1\pi_{12}(p)}\subset \R_T^2\rightarrow \overline{\alpha_1\pi_{31}(p)}\subset \R_T^2$ is a bi-Lipschitz bijection, with constant independent of $\eps$. A direct computation then provides that the map $\Psi:\overline{\zeta^1w^a}\subset \R_S^2\rightarrow \overline{\zeta^1w^c}\subset\R_S^2$ defined by
\begin{align}\label{omeo}
 \Psi:=(L_{31})^{-1}\circ \pi_{31}\circ \Phi_{12}\circ L_{12},
\end{align}
is  bi-Lipschitz between $\overline{\zeta^1w^a}$ and $\overline{\zeta^1w^c}$ with bi-Lipschitz  constants of order $1$ as $\eps\to 0^+$. Given $Q\in \overline{\zeta^1w^a}$ let  $Q':=\Psi(Q)\in \overline{\zeta^1w^c}$.

Now we show that $T_1^\eps$ is 
foliated by the segments $\overline{QQ'}$, {\it i.e.}, for any $ R\in T_1^\eps$ we can find a unique $Q \in \overline{\zeta^1w^a}$ for which $R \in \overline{QQ'} $. 

First we notice that 
$\overline{QQ'}\cap \overline{SS'}=\emptyset $ for any $ Q\not = S\in \overline{\zeta^1w^a}$ with $ Q'=\Psi(Q)$ and $ S'=\Psi(S).$
Indeed, thanks to the fact that $\Psi$ is a homeomorphism and that it keeps $\zeta^1$ fixed, it is 
easy to see that if $S\in \overline{\zeta^1Q}$ then $S'\in \overline{\zeta^1Q'}$, or if $Q\in \overline{\zeta^1S}$ then $Q'\in \overline{\zeta^1S'}$. Consider the function 
\begin{equation*}
f(q,\sigma)=q\tau+\sigma \nu(q),\qquad q\in [0,|w^a-\zetai|], \sigma \in [0,|\Psi(q\tau)-q\tau|],
\end{equation*}
where $\tau:=\frac{w^a-\zetai}{|w^a-\zetai|}$ and $\nu(q):=\frac{\Psi(q\tau)-q\tau}{|\Psi(q\tau)-q\tau|}$. It is clear that the image of $f$ is a closed set and ${\rm Im}(f)=\{\overline{QQ'}:Q\in \overline{\zeta^1w^a},Q'=\Psi(Q)\}.$ Now we show that
$
{ \rm Im}(f)=T^\eps_1.
$
Assume by contradiction there is $R\in T^\eps_1\setminus  {\rm Im}(f)$ and take a disk $B\subset T^\eps_1\setminus { \rm Im}(f)$ centered at $R$. Let $Q_r, Q_l\in \overline{\zeta^1w^a}$ be such that $q_r:=|Q_r-\zetai|$ (resp. $q_l=|Q_l-\zetai|$) be the supremum (resp. the infimum) parameter for which $B$ lies on the right (resp. left) of $\overline{Q_r Q_r'}$ (resp. $\overline{Q_l Q_l'}$). Note that $Q_r \neq Q_l$ due to the injectivity of $\Psi$, thus for any $Q\in \overline{Q_r Q_l}$ the segment $\overline{Q Q'}$ must intersect $B$, a contradiction, see 
Figure \ref{fig:Ballright}. 
\begin{figure}
\centering
\begin{subfigure}{.5\textwidth}
  \centering
  \begin{center}
           \includegraphics[scale=0.35]{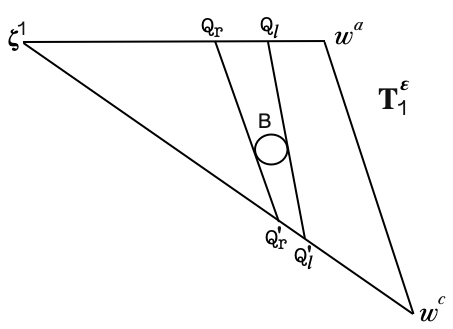}  
  \end{center}
 \caption{}
  \label{fig:Ballright}
\end{subfigure}%
\begin{subfigure}{.5\textwidth}
  \centering
  \begin{center}
   \includegraphics[scale=0.35]{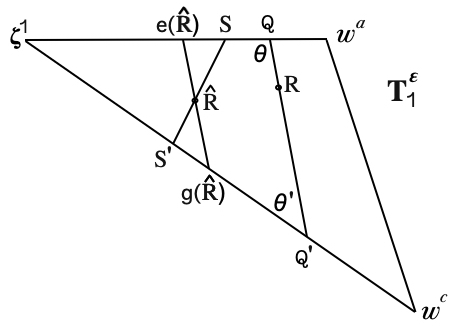}  
      \end{center}
  \caption{}
  \label{fig:thetapic}
\end{subfigure}
\caption{Surjectivity of the foliation in $T^\eps \subset \R^2_S$.}
\end{figure}
Hence we may define $\uh$ on $T_1^\eps$ as
\begin{align}\label{cavolo}
\uh(R):=\uh(Q),\qquad R\in \overline{QQ'}, ~Q\in \overline{\zeta^1w^a}.
\end{align}
We want now to show that on $T_1^\eps$, $\uh$ is Lipschitz continuous with Lipschitz constant of order $\frac{1}{\eps}$. To prove this let us fix $R\in T^\eps_1$. By definition $\uh(R)=\uh(Q)$ for some $Q\in \overline{\zeta^1w^a}$ and $\uh$ is constant on the segment $\overline{QQ'}\ni R$.

Let $\e: T^\eps_1 \to \overline{\eti w^a}$ be the function taking $(x,y)\in T^\eps_1$ to the intersection point of $\overline{\eti w^a}$ and the line passing through $(x,y)$ parallel to $\overline{QQ'}$. 
Let $\g: T^\eps_1\to \overline{\eti w^c}$ be the function taking $(x,y)\in T^\eps_1$ to the intersection point of $\overline{\eti w^c}$ and the line passing through $(x,y)$ parallel to $\overline{QQ'}$. Let $\hat R\in T^\eps_1$ a point in $T^\eps_1$, we want to estimate the ratio
$$\frac{\vert \uh(\hat R)-\uh(R)\vert}{|\hat R- R|}.$$
Consider the two segments $\overline{Q \e(\hat R)}$ and $\overline{Q\g(\hat R)}$. By definition $\hat R\in \overline{SS'}$ and $\uh(\hat R)=\uh(S)=\uh(S')$ for two points $S\in \overline{\eti w^a}$ and $S'\in \overline{\eti w^c}$. It is straightforward that either $S\in \overline{Q \e(\hat R)}$ or $S'\in \overline{Q'\g(\hat R)}$. Without loss of generality suppose the first case holds, see Figure \ref{fig:thetapic}. 

Finally, denote by $\theta$ the angle between $\overline{QQ'}$ and $\overline{\eti w^a}$ and by $\theta'$  the angle between $\overline{QQ'}$ and $\overline{\eti w^c}$. Using the fact that the homeomorphism in \eqref{omeo} is bi-Lipschitz with constant of order $1$ it is not difficult to see that there is a constant $\theta_0>0$ independent of $\eps$ such that $\min\{\theta,\theta'\}\geq\theta_0$. 
This is a consequence of the fact that the bi-Lipschitz constant of $\Psi$ in \eqref{omeo} is of order $1$. Indeed, if $L={\rm lip}(\Psi)$ and $1/L'={\rm lip}(\Psi^{-1})$, we see that $$L' \leq\frac{|Q'-\zeta^1|}{|Q-\zeta^1|}\leq L,$$ hence $$\frac{1/L+\cos \theta_{\zeta^1}}{\sin{\theta_{\zeta^1}}}\leq \frac{\cos \theta}{\sin \theta}\leq \frac{1/L'+\cos \theta_{\zeta^1}}{\sin{\theta_{\zeta^1}}},$$
where $\theta_{\zeta^1}$ is the angle at $\zeta^1$ (here we have used the law of sines and that $\theta'=\pi-\theta_{\zeta^1}-\theta$). A similar estimate holding for $\theta'$, this readily provides the boundedness from below of $\min\{\theta,\theta'\}$.

As a consequence we have  
\begin{align*}
 |\hat R- R|\geq|Q-\e(\hat R)|\vert\sin\theta\vert\geq |Q-\e(\hat R)|\vert \sin\theta_0\vert. 
\end{align*}
Thus, we compute
\begin{align}
 \frac{\vert\uh(\hat R)-\uh(R)\vert}{|\hat R- R|}\leq\frac{\vert\uh(Q)-\uh(S)\vert}{|Q-\e(\hat R)|\vert\sin\theta\vert }\leq \frac{\vert\uh(Q)-\uh(S)\vert}{|Q-S|\vert\sin\theta\vert }\leq \frac{1}{\vert\sin\theta_0\vert}\frac{\vert\uh(Q)-\uh(S)\vert}{|Q-S| },
\end{align}
that is bounded by the Lipschitz constant of $\Phi_{12}\circ L_{12}$ which is of order $\frac{1}{\eps}$.

Eventually we compute the Jacobian of $\uh$ in \eqref{cavolo}. By construction the image of $T^\eps_1$ by $\uh$ is exactly the curve $\Gamma_1$, which has zero Lebesgue measure in $\R^2$. By a standard application of the area formula it follows that the Jacobian of $\uh$ is vanishes a.e. in $ T^\eps_1$. We have concluded the definition of $\uh$ in $T^\eps_1$.
The constructions on $T^\eps_2$ and on $T^\eps_3$ are similar, and similar estimates of the derivatives and Jacobian hold. 

{Using that the area of the triangle $T^\eps$ is of order $\eps^2$, we have}
\begin{equation}\label{Teps_contribution}
\lim_{\eps \to 0} 
\A(\uh, T_\eps) = \lim_{\eps \to 0} O(\eps)+O(\eps^2)=0.
\end{equation}

From \eqref{def:uhregioni}, \eqref{def:uhstrip},
\eqref{T0_contribution}, \eqref{cavolo}, and the estimates above
it follows that 
\begin{equation}\label{ultima}
\{\uh\} \subset {\rm Lip}(\sourcedisk; \R^2), \qquad \qquad
\lim_{\eps\to 0} \int_{\sourcedisk} \vert \uh - u\vert~dxdy=0.
\end{equation}
Moreover
\begin{equation}
\A(\uh, \sourcedisk) = \A(\uh, \Aae\cup \Bbe \cup \Cce)
+ 
\A(\uh, S^\eps_\ac) + 
\A(\uh, S^\eps_\bc) + 
\A(\uh, S^\eps_\cc) + 
\A(\uh, T^\eps).\label{wecansplit}
\end{equation}
Then \eqref{eqn:case1f} follows from  \eqref{wecansplit}, \eqref{areanulla}, \eqref{3cygar_contribution}, \eqref{chegabbarmivolete} and \eqref{Teps_contribution}.

\medskip

{\it Case 2.} Assume that two of the segments separating $\Aa,~ \Bb,~ \Cc$ meet at the origin with an angle greater than or equal to $\pi $.

Similar to Case 1, we divide the domain $\sourcedisk$ into a finite number of subsets and define the sequence $\{\uh\}$ on each of these sets. Draw the normal to each segment at the point of distance $\del$ from the origin. The normal lines meet at two points $\eti$, $\etii$. Divide $\sourcedisk$ 
into three cygar-shape subsets $S_\ac^\eps,~S_\bc^\eps,~S_\cc^\eps$ (with widths of order $\delta=O(\eps)$) and a quadrilateral $H^\eps$ as in Figure \ref{fig:D180}. 
\begin{figure}
\centering
\begin{subfigure}{.5\textwidth}
  \centering
  \begin{center}
           \includegraphics[scale=0.45]{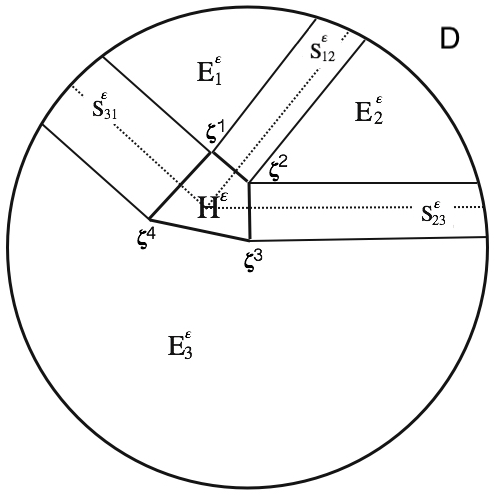}  
  \end{center}
 \caption{The non acute case.}
  \label{fig:D180}
\end{subfigure}%
\begin{subfigure}{.5\textwidth}
  \centering
  \begin{center}
   \includegraphics[scale=0.4]{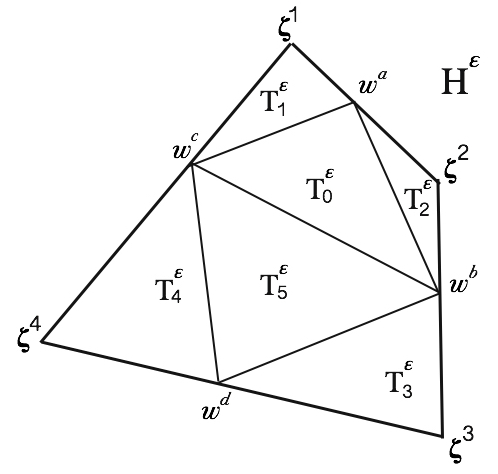}  
      \end{center}
  \caption{Zoom of $H^\eps$ in (a).}
  \label{fig:H_eps}
\end{subfigure}
\caption{}
\end{figure}
Let
$$\Aae
:=
\Aa \setminus \left(S_\bc^\eps \cup H^\eps \cup S_\cc^\eps \right),
\quad
\Bbe := 
\Bb \setminus \left(S_\ac^\eps \cup H^\eps \cup S_\cc^\eps \right), 
\quad \Cce := 
\Cc \setminus \left(S_\ac^\eps \cup H^\eps \cup S_\bc^\eps \right).
$$
Set
\begin{equation}
\uh := 
\begin{cases}
\alp & {\rm in}~ \Aae,
\\
\bet & {\rm in}~ \Bbe,
\\
\gam & {\rm in}~ \Cce.
\end{cases}
\end{equation}
Define $\uh$ on $ S_\ac ^\eps \cup S_\bc ^\eps \cup S_\cc ^\eps$ as in 
Step 2 in case 1. It remains to define $\uh$ on $H^\eps$.
Recall that by construction there exist uniquely determined three points $w^a \in \overline{\zeta^1\zeta^2} $, $w^b \in\overline{\zeta^2\zeta^3}$ and $w^c\in \overline{\zeta^1\zeta^4}$ such that $$\uh(w^a)=\uh(w^b)=\uh(w^c)=\targetorigin.$$
Divide $H^\eps$  into six triangles $T^\eps_0,~T^\eps_1,~T^\eps_2,~T^\eps_3,~T^\eps_4,~T^\eps_5$, as in Figure \ref{fig:H_eps}, where $w^d$ is any point in $\overline{\zeta^3 \zeta^4}$ and $w^d\not =\zeta^3,~ w^d\not =\zeta^4$.

Set 
$$\uh:=\targetorigin \qquad {\rm in} ~~~T^\eps_0 .$$
We define $\uh$ in the  triangles $T^\eps_1$ and $T^\eps_2$ as in 
Step 3; it remains to define $\uh$ on $T^\eps_3,~T^\eps_4,~T^\eps_5$.
Let us first define $\uh$ on the edges $\overline{w^cw^d}$ and $\overline{w^bw^d}$. The map $\uh$ is already defined on the other edges, and its graph over $\overline{\zeta^4w^c}$ and $\overline{\zeta^3w^b}$ is given by a suitable reparametrization of the curve $\Gamma_3$, whereas $\uh$ on $\overline{\zeta^4\zeta^3}$ is constantly $\alpha_3$. Therefore it suffices to define $\uh$ in such a way its graph over $  \overline{w^cw^d}$ and $\overline{w^bw^d}$ coincides with $\Gamma_3$ as well, and then we can define $\uh$ inside  $T^\eps_3$ and $T^\eps_4$ using the same construction for $T^\eps_1$ in 
step 3. Similarly, using that the graph of $\uh$ on  $  \overline{w^cw^d}$ and $\overline{w^bw^d}$ is again $\Gamma_3$, we can repeat the construction in the triangle $T^\eps_5$. 
Following the computation as in case 1 we get \eqref{eqn:case1f}.
This concludes the proof. 

\end{proof}


\begin{proof}[Proof of Theorem \ref{thm:upper_bound_with_infimum}]
We will suitably adapt the construction made in the proof of Proposition \ref{prop:piecewise_linear}. By hypothesis  the regions $\Aa, \Bb, \Cc$ are enclosed by $C^2$-embedded curves $C_{ij},\;ij\in\{12,23,31\},$ parametrized by arc length $\scij:[0,r_{ij}]\to \R^2 ,~ ij\in \{\ac,\bc,\cc\}$. Moreover such curves meet $\partial\sourcedisk$ transversely and intersect each other (transversely) only at one point $\neworigin$. Suppose that the angles formed at $\neworigin$ by the three curves are all less than $\pi$ (the other case is similarly adapted from the corresponding case in the proof of Proposition \ref{prop:piecewise_linear}). 
We will divide the domain $\sourcedisk$ into a finite number of subsets and define the sequence $\{\uh\}$ on each of these sets.

Let $\delta_\eps>0$ be such that $\del_\eps \to 0 $ as $\eps \to 0^+$. Let $\z\in[0,r_{ij}]$ be an arc lenght parameter on $\cij$, with orthogonal coordinate $d$ that coincides with the signed distance from $\cij$ negative in $E_i$ and positive in $E_j$. Let $Q_{ij}\in \cij$ be the point with arc distance $\z=\delta_\eps$ from the origin $\neworigin$. Consider the three lines normal to $\cij$ at $Q_{ij}$. For $\delta_\eps$ sufficiently small, since the angles at the origin are less than $\pi$ and the curves are of class $C^2$ up to the closure, these lines mutually meet at points $\zeta^1$, $\zeta^2$, and $\zeta^3$. Let $\eps_{ij}$ be the length of $\overline{\zeta^i\zeta^j}$, which are of order ${\eps}$. The tubular coordinates of the points $\zeta^1$ and $\zeta^2$ with respect to $C_{12}$ are $(d_1,\delta_\eps,)$ and $(d_2,\delta_\eps)$, with $d_2-d_1=\eps_{12}, d_1<0, d_2>0$. For $\delta_\eps$ small enough we can consider the cylindrical neighborhood of $C_{12}$ defined as
\begin{align}
 S_{12}^\eps:=\{(x,y)\in D:\z(x,y)\geq \delta_\eps,\;d(x,y)\in (d_1,d_2)\},
\end{align}
where we have prolonged $C_{12}$ outside $D$ for convenience.
Similarly we define $S_{23}^\eps$ and $S^\eps_{31}$. Let $T^\eps $ be the triangle with vertices $\zeta^1$, $\zeta^2$, and $\zeta^3$.

\begin{figure}
\centering
\begin{subfigure}{.7\textwidth}
  \centering
  \begin{center}
           \includegraphics[scale=0.35]{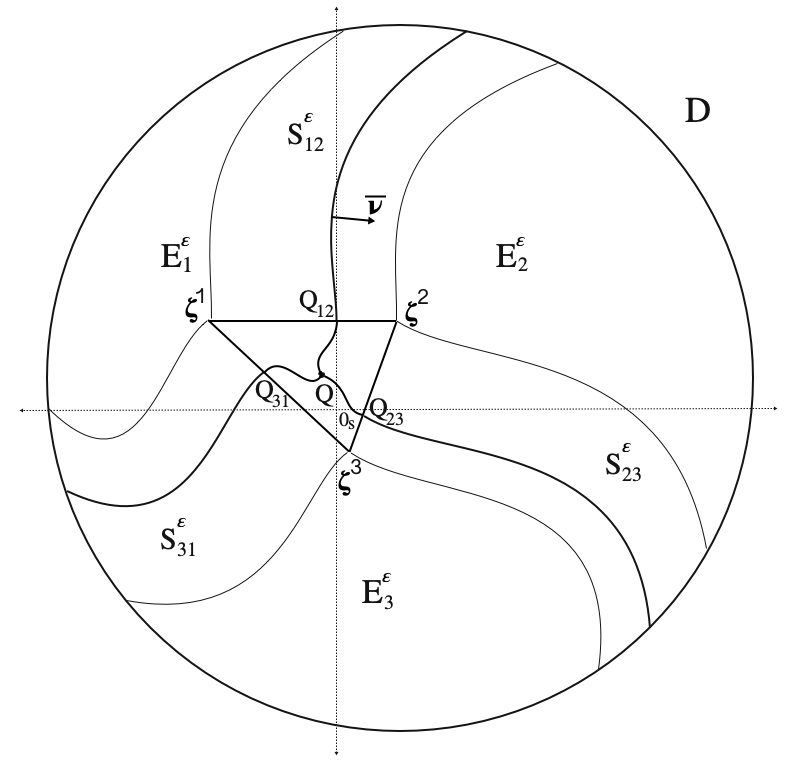}  
  \end{center}
 \caption{}
  \label{fig:Dcurves2}
\end{subfigure}%
\begin{subfigure}{.3\textwidth}
  \centering
  \begin{center}
   \includegraphics[scale=0.5]{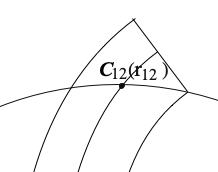}  
      \end{center}
  \caption{}
  \label{fig:roc}
\end{subfigure}
\caption{proof of Theorem \ref{thm:upper_bound_with_infimum}.}
\end{figure}

Finally, let $\Aae, ~\Bbe, ~\Cce$ be defined as in \eqref{eqn:4.51}, and $\uh$ as in \eqref{def:uhregioni}.
 
\smallskip
{\it Step 1.} Definition of $\uh$ on $S^\eps_\cc\cup S^\eps_\ac \cup S^\eps_\bc$. 
We do the construction on $S^\eps_\cc$, and $\uh$ will be defined similarly on $S^\eps_\ac$ and $ S^\eps_\bc$. We know that $\sGammacc ([\delta_\eps,r_\cc])=\Gammacc \cap S^{\eps}_\cc$. The system of  coordinates $(d,\z)$ defines a $C^1$-diffeomorphism $h$ between the rectangle $[d_1,d_2]\times[\delta_\eps,\rho_\cc]$ and its image $\NNiii$ which contains $S^\eps_{12}$, namely 
\begin{align*}
\h:[d_1,d_2]\times[\delta_\eps,\rho_\cc] \to \NNiii;\qquad
\h(\w,\z):=\sGammacc(\z)+\w \nub(\z),
\end{align*}
where $\nub(\z)$ is the unit normal vector pointing toward $E_2$ at $\sGammacc(\z)$ and $\ro_\cc=\ro_\cc^\eps \geq r_\cc$ is the infimum of those $\ro$ for which $S^\eps_\cc\subset\NNiii$, see Figure \ref{fig:roc}. 

Since $\h$ is a $C^1$-diffeomorphism we have that 
\begin{align*}
\hinv:\NNiii\to [d_1,d_2]\times[\delta_\eps,\rho_\cc];\qquad
\hinv(x,y):=(\w(x,y),\z(x,y)),
\end{align*}
is the inverse of $\h$ and is of class $C^1$. We want to estimate the Jacobian of $h^{-1}$. To this aim, we first see that $\nabla \w ({\sGammacc(\z)})=\nub(\z)$ since $\sGammacc([\delta_\eps,\ro_\cc])$ is the zero level set of $\w$ and, from \cite[Rem. 3(1)]{Amb:00}, we have
\begin{equation}\label{grad d}
\nabla \w ({\h(\w,\z)})= \nub(\z),\qquad (\w,\z)\in [d_1,d_2]\times[\delta_\eps,\rho_\cc].
\end{equation}
Fix $\z \in [\delta_\eps,\rho_\cc]$; by definition of tubular coordinates the segment $\{\sGammacc(\z)+\w \nub(\z):\w \in [\w_1,\w_2]\}$ is a level set of the function $\z(\cdot)$, hence, 
\begin{equation}
\nabla\z(x,y) \perp \nub(\z(x,y)), \label{eqn:4.271}
\end{equation}
therefore 
\begin{equation}\label{grad d grad z}
\nabla \z ~{\Huge \cdot}~ \nabla \w~ =~0\qquad \text{ in }S^{\eps}_{12}.
\end{equation} 
Thus the Jacobian of $h^{-1}$ will be
\begin{align}\label{jac}
j(h^{-1})=|\nabla \z||\nabla d|=|\nabla \z|, 
\end{align}
since $|\nabla d|=1$ in $\NNiii$.
Let us compute $\nabla \z$; fix $\w \in (d_1,d_2)$ and define $\sGammacc^{\w} (\z):=\sGammacc(\z)+\w \nub(\z)$. 
Now recall \eqref{eqn:4.271} and that
$(\sGammacc^{\w})' (\z)$ is parallel to $\bar \nu^\perp(\z)$, so that 
\begin{equation}\label{gradz1}
|\nabla \z|=\nabla \z\cdot\bar\nu^\perp=\frac{\nabla \z\cdot (\sGammacc^{\w})'}{|(\sGammacc^{\w})'|}\qquad \text{ in }S^{\eps}_{12}. 
\end{equation}
Let us recall that $\Gammacc$ is parametrized by arc length, {\it i.e.}, $|\sGammacc^\prime(\z)|=1$, so that  $\nub^{\prime}(\z)=|\sGammacc^{\prime\prime}(\z)|\sGammacc^\prime(\z)$. Thus ${(\sGammacc^\w)}^{\prime}(\z)=\left(1+\w |\sGammacc^{\prime\prime}(\z)|\right)\sGammacc^{\prime}(\z)$. Since $\z\circ \sGammacc^\w={\rm Id}$ it follows that $\nabla \z (\sGammacc^\w(\z))^T (\sGammacc^\w)^\prime(\z)=\nabla \z (\sGammacc^\w(\z))\cdot (\sGammacc^\w)^\prime(\z)=1$. Therefore, from \eqref{gradz1}, we deduce 
\begin{equation}\label{grad z}
|\nabla \z|= \frac{1}{1+\w|\sGammacc^{\prime\prime}(\z)|}\qquad \text{ in }S^{\eps}_{12},
\end{equation} 
and in particular $\lim_{\w\rightarrow 0} |\nabla \z| =1$ uniformly in $S^{\eps}_{12}$. 

We are ready to define $\uh$ in $S^\eps_{12}$. We first set $\psi_{12}^\eps$ as in \eqref{kappaeps} with $r_{12}+c_\eps=\rho_{12}$, {\it i.e.}, $\psi_{12}^\eps(\z)=\keps(z-\delta_\eps),$ setting $\keps:=\frac{\rho_{12}}{\rho_{12}-\delta_\eps}$. Then we  define $\tilde u ^\eps$ on $[d_1,d_2]\times [\delta_\eps,\rho_{12}]$ as in the right hand side of \eqref{def:uhstrip} and set
\begin{align}
 \uh:=\tilde u ^\eps\circ h^{-1}\qquad \text{ in }S^{\eps}_{12}.
\end{align}
Explicitly, recalling that $\xi=\frac{\alpha_2-\alpha_1}{\ell_{12}}$ and $\eta=\xi^\perp$, for $ (x,y) \in S^\eps_\cc$ we have
\begin{equation}\label{uhscurves}
\uhs(x,y) := \alp + \left(\frac{\w(x,y)-d_1}{\eps_\cc} \right)\ell_\cc  \xi + \tenda_\cc^{\sigma}\left(\frac{\w(x,y)-d_1}{\eps_\cc}\ell_{12}
~,~\keps\left(\z(x,y)-\delta_\eps\right)\right) \eta.
\end{equation}
Observe that $\uhs = (\uhs_1, \uhs_2) \in \Lip(S^\eps_\cc; \R^2)$,
$\uhs =\alp$ on $\{(x,y) \in S^\eps_\cc : \w(x,y)=d_1\}$, $\uhs = \bet$ on $\{(x,y) \in S^\eps_\cc : \w(x,y) = d_2\}$, and by construction there exists $w_a\in \h\left([d_1,d_2]\times\{\delta_\eps\}\right)$ such that $\uhs(w_a)=0_{T}$.
Write for simplicity $\widetilde \tenda = 
\tenda^{\sigma}_\cc.
$
We have		
\begin{align*}
\grad \uhs_1 = \left(\frac{\ell_\cc \xi^1}{\eps_\cc}\w_x 
+\frac{\ell_\cc \eta^1}{\eps_\cc}\widetilde \tenda_s \w_x +\keps\widetilde \tenda_t \z_x~,~ 
\frac{\ell_\cc \xi^1}{\eps_\cc}\w_y 
+ \frac{\ell_\cc \eta^1}{\eps_\cc}\widetilde \tenda_s \w_y +\keps\widetilde \tenda_t \z_y
\right),\\ 
\grad \uhs_2 
= \left(\frac{\ell_\cc \xi^2}{\eps_\cc}\w_x 
+ \frac{\ell_\cc \eta^2}{\eps_\cc}\widetilde \tenda_s \w_x +\keps\widetilde \tenda_t \z_x~,~ 
\frac{\ell_\cc \xi^2}{\eps_\cc}\w_y 
+ \frac{\ell_\cc \eta^2}{\eps_\cc}\widetilde \tenda_s \w_y +\keps\widetilde \tenda_t \z_y
\right),
\end{align*}
where 
$\widetilde \tenda_s, \widetilde \tenda_t$ denote the partial derivatives
of $\widetilde \tenda$ with respect to $s=\frac{\w(x,y)-d_1}{\eps_\cc}\ell_{12}$ and $t=\keps(\z(x,y)-\delta_\eps)$ respectively, and 
are evaluated
at $\left(\frac{\w(x,y)-d_1}{\eps_\cc}\ell_{12}~,~ \keps(\z(x,y)-\delta_\eps)\right)$. 
Hence
\begin{equation*}
\vert \grad \uhs_1\vert^2 + \vert
\grad \uhs_2\vert^2
~ =~    
\frac{\ell_\cc^2}{\eps_\cc^2}|\nabla \w|^2+
\frac{\ell_\cc^2}{\eps_\cc^2}|\nabla \w|^2 (\widetilde \tenda_s)^2+\keps^2|\nabla \z|^2 (\widetilde \tenda_t)^2+\frac{2\ell_\cc }{\eps_\cc}\keps\left(  \nabla \w \cdot \nabla \z \right) \widetilde \tenda_s \widetilde \tenda_t
\end{equation*}
where 
 we have used $\vert \xi\vert= \vert \eta\vert =1$ and $\xi_1 \eta_1 + 
\xi_2 \eta_2=0$. From \eqref{grad d grad z} we have
\begin{equation} \label{contocruccurves}
\vert \grad \uhs_1\vert^2 + \vert
\grad \uhs_2\vert^2
~ =~    
\frac{\ell_\cc^2}{\eps_\cc^2}+\frac{\ell_\cc^2}{\eps_\cc^2} (\widetilde \tenda_s)^2+\keps^2|\nabla \z|^2 (\widetilde \tenda_t)^2.
\end{equation}
Moreover
\begin{equation}\label{determincurves}
\begin{aligned}
\left(\frac{\partial \uhs_1}{\partial x}
\frac{\partial \uhs_2}{\partial y}-\frac{\partial \uhs_1}{\partial y}
\frac{\partial \uhs_2}{\partial x}\right)^2
 ~=~ & 
\frac{\ell_\cc^2 }{\eps_\cc^2}\keps^2 
(\widetilde \tenda_t)^2
 \left(\w_x\z_y-\w_y\z_x\right)^2 (\xi^1\eta^2 - \xi^2 \eta^1)^2= 
\frac{\ell_\cc^2 }{\eps_\cc^2}\keps^2 
(\widetilde \tenda_t)^2
 |\nabla \z|^2  ,
\end{aligned}
\end{equation}
where again $\widetilde \tenda_s, \widetilde \tenda_t$ are evaluated
at $\left(\frac{\w(x,y)-d_1}{\eps_\cc}\ell_{12}~,~ \keps(\z(x,y)-\delta_\eps)\right)$, and 
we have used \eqref{grad d grad z}, \eqref{jac}, and
$\xi^1\eta^2 - \xi^2 \eta^1= 1$.
Therefore from \eqref{contocruccurves} and \eqref{determincurves} we obtain
$$
\begin{aligned}
& 1 + \vert \grad \uhs_1\vert^2 + \vert
\grad \uhs_2\vert^2 + \left(\frac{\partial \uhs_1}{\partial x}
\frac{\partial \uhs_2}{\partial y}-\frac{\partial \uhs_1}{\partial y}
\frac{\partial \uhs_2}{\partial x}\right)^2
\\
= & 
1 + 
\frac{\ell^2_\cc}{\eps_\cc^2} \left( 
1+ \big(
\widetilde \tenda_s \big)^2 
+ 
\big(
\widetilde \tenda_t
\big)^2
\keps^2 \left(1+
\frac{\eps_{12}^2}{\ell_{12}^2}
\right)|\nabla \z|^2
\right).
\end{aligned}
$$
As a consequence
\begin{equation}\label{contributonellastrisciacurves}
\begin{aligned}
&\area(\uh, S^\eps_\cc)
 =~  \frac{\ell_\cc}{\eps_\cc}
 \int_{S_\cc^\eps}
\sqrt{
1 + \left(
\widetilde \tenda_s
\right)^2 
+ 
\left(\widetilde \tenda_t\right)^2
\keps^2\left(1+\frac{\eps_{12}^2}{\ell_{12}^2}\right)
|\nabla \z|^2 
+ O(\eps^2)
}~dxdy\\
=&  
\frac{1}{\keps} \int_{\rettangolo_\cc\setminus P_\eps}
\frac{1}{|\nabla \z|} 
\sqrt{
1 + 
\left(\widetilde \tenda_s\left(s,t\right)\right)^2 
+ 
\left(\widetilde \tenda_t\left(s,t\right)\right)^2
\keps^2 
\left(1+\frac{\eps_{12}^2}{\ell_{12}^2}\right)
|\nabla \z|^2 
 + O(\eps^2)
}
~dsdt,
\end{aligned}
\end{equation}
where $\widetilde \tenda_s, \widetilde\tenda_t$ in the first integral are evaluated at $ \left(
 \frac{\w(x,y)-d_1}{\eps_\cc}\ell_{12}~,~ \keps(\z(x,y)-\delta_\eps) 
\right)$, $\nabla \z$ in the second integral is evaluated at $(x,y)=\Phi^{-1}(s,t)$ and  the last equality follows from the change of variables
$$
\Phi: (x,y) 
\in \NN_\cc^{\eps} \to 
\left( \frac{\w(x,y)-d_1}{\eps_\cc}\ell_{12}~,~\keps(\z(x,y)-\delta_\eps) 
\right) 
= (s,t)
\in \rettangolo_\cc.$$ 
and $P_\eps := \rettangolo_\cc \setminus \Phi(S^\eps_\cc) $
(see Figure \ref{fig:Peps}). Here one checks that $\Phi=H\circ h^{-1} $ with $H(d,\z)= (\frac{\w-d_1}{\eps_\cc}\ell_{12}~,~ \keps(\z-\delta_\eps))$ so that, using \eqref{jac}, the Jacobian of the change of variable is $\frac{1}{|\nabla \z(\Phi^{-1}(s,t))|}\frac{\eps_{12}}{\ell_{12}\keps}$. 
Hence, recalling \eqref{grad z} and that $\keps \to 1$ as $\eps\to0^+$, 
\begin{equation}\label{contributonellastriscialimitecurves}
\lim_{\eps \to 0^+} 
\area(\uhs, S^\eps_\cc) = 
 \int_{\rettangolo_\cc} 
\sqrt{
1 + 
\big(\widetilde \tenda_s
\big)^2
+ 
\big(\widetilde \tenda_t
\big)^2
}
~dsdt.
\end{equation}
Now, let us recall that $\tilde m=m_{12}^\sigma$ is the approximating function as in \eqref{lip approx m}; it follows that
\begin{equation}\label{isitimmediatecurves}
 \int_{\rettangolo_\cc} 
\sqrt{
1 + 
\big(\widetilde \tenda_s
\big)^2
+ 
\big(\widetilde \tenda_t
\big)^2
}
~dsdt = 
{\mathfrak A}_\cc(\Gamma) + O(\sigma).
\end{equation}
Hence, employing the same construction in the strips $S^\eps_\ac$  and $S^\eps_\bc$, 
and using \eqref{isitimmediatecurves} we obtain from a diagonal argument with $\sigma=\sigma_\eps \rightarrow0$ as $\eps\rightarrow 0^+$,
\begin{equation}\label{cruccurves}
\lim_{\eps \to 0^+} 
\area(\uhs, 
S^\eps_\ac \cup 
S^\eps_\bc \cup 
S^\eps_\ac 
) = {\mathfrak A}_{12}(\Gamma)+{\mathfrak A}_{23}(\Gamma)+{\mathfrak A}_{31}(\Gamma).
\end{equation}
\smallskip

{\it Step 2.} Definition of $\uh$ on $ T^\eps$. This is identical to 
Step 3 of the proof of Proposition \ref{prop:piecewise_linear} and therefore $\{\uh\} \subset {\rm Lip}(B_\raggio; \R^2)$ and \eqref{ultima} holds. Following the same computations of Proposition \ref{prop:piecewise_linear} the conclusion follows.

\smallskip
{\it Step 3.}
For the case where two of the curves $\cij, ij\in \{12,23,31\}$ meet at $Q$ with an angle larger than or equal to $\pi$ we replace $T^\eps$ with $H^\eps$ defined in case 2 of Proposition \ref{prop:piecewise_linear}, in the above construction. 

\end{proof}
%

\section{Existence of minimizers for the functional $\AG$}\label{sec:AnalysisAG}
Let $\sourcedisk$ be an open disk centered at the origin such that $\Aa,~\Bb,~\Cc$ are circular sectors with $120^\circ$ angles and let $\targettriangle$ be an equilateral triangle. Let $\genericp$ be the barycenter of $\targettriangle$ and $\pbarGammaalpi$ be the segment connecting $\alpi$ and $\genericp$, $i\in\{1,2,3\}$. Hence $\barelementXclip=(\barelementXclip_1,\barelementXclip_2,\barelementXclip_3)\in X_\Lip$ so that 
\begin{equation*} \label{eqn:bel}
\inf\big \{\AG(\elementXclip): \elementXclip\in X_\Lip\big \} \leq \AG(\barelementXclip).
\end{equation*}
Moreover we have
\begin{equation*} \label{eqn:scala}
\vert \sourcedisk \vert + \AG(\barelementXclip)=\A(u,\sourcedisk)\leq \vert \sourcedisk \vert+\inf\big \{\AG(\elementXclip): \elementXclip\in X_\Lip\big \},
\end{equation*} 
where $u=u_{\rm symm}$ (see Section \ref{sec:introduction}), and the equality follows from \cite[Section 3]{Sc:19} 
and the inequality follows from Proposition \ref{prop:piecewise_linear}. Thus 
$$\AG(\barelementXclip)= \min\big \{\AG(\elementXclip): \elementXclip\in X_\Lip\big \}.$$
Hence in this symmetric situation the optimal connection is obtained through the Steiner graph connecting $\alp,$ $\bet$ and $\gam$. This motivates the analysis of this section, which is carried on without symmetry assumptions.

{We recall that given a connection $\Gamma=(\Gamma_{1},\Gamma_{2},\Gamma_{3})\in X$ we denote by $\varphi_{ij}=\varphi_{ij}(\Gamma_{ij}):[0,\ellalpij]\rightarrow \R$ the function whose graph is $\Gamma_{ij}=\Gamma_{i}\cup\Gamma_{j}$ (see \eqref{varphi_ij}). }
\begin{Definition}[\textbf{Convergence in $X$}]\label{ConvinX}
We say that a sequence $\{\elementXclip^n\} \subset X$ converges to $\elementXclip \in X$ in $X$, and we write $\Gamma^n\rightarrow\Gamma$ in $X$,  if
\begin{equation}
\phiij(\pnGammaalpij)\to \phiij(\pGammaalpij) {~~~\rm in~}\Lone([0,\ellalpij]), \qquad ij \in \{12,23,31\}.
\end{equation}
\label{Def:XConv}
\end{Definition}
\subsection{Density and approximation}
We start to show that a $\BV$ connection $\elementXclip \in X$ can be approximated by Lipschitz connections; the difficulty is to keep graphicality of each branch of of the approximating connections
with respect to the two corresponding edges of $\targettriangle$ at the same time.

Recall that $\Gamma_i$ is the branch of the connection $\Gamma$ connecting $\alpha_i$ to $p$ and that by Definition \ref{def:connections} we have
$$ 
\Gamma_i\res \targettriangle \setminus \overline{p \pi_{ij}(p)} \cup \overline{p \pi_{ki}(p)} ={\Gamma}_{\varphi_{ij}\res[0,w_{ij})}={\Gamma}_{\varphi_{ki}\res(w_{ki},\ell_{ki}]}.
$$
Note that we excluded the vertical parts over the points $\pi_{ij}(p),~ij\in\{12,23,31\},$ due to Remark \ref{Rmk:graph.}; however we still have 
$$
\Gamma_i \cup \Gamma_j={\Gamma}_{\varphi_{ij}}.
$$

\begin{Lemma}[\textbf{Piecewise linear approximation}]\label{Lem:density}
For any $\elementXclip \in X$ with target triple point $\genericp \in \targettriangle$ there exists a sequence $\{\elementXclip^n\}\subset \Xclip$ of connections with target triple point $\genericp$ such that $\phiij(\pnGammaalpij), ij\in \{12,23,31\},$ is a 
piecewise linear\footnote{This means that it is Lipschitz piecewise 
linear with at most finitely many points of nondifferentiability.} 
function, 
$$
\Hone(\pnGammaalpij)\leq\Hone(\Gammaalpij)
$$
 and 
\begin{equation}
\elementXclip^n \to\elementXclip\qquad {\rm ~in~} ~X. 
\label{eqn:densityconv}
\end{equation}
\end{Lemma}
\begin{proof} Let $ij=12$ and let $w_{12}$ be defined as in \eqref{eqn:wij}. 
Let $\nalpbet:=(0,1)\in\R^2$ be the inward unit normal to $\overline{\alp \bet}$, $\ngamalp:=(\alpha,\beta)$ be the inward unit normal to $\overline{\gam \alp }$, and $\nualp(\bar s):=(\nualp^1(\bar s), \nualp^2(\bar s))$ be the generalized outward unit normal at the point $(\bar s,\varphi_{12}(\bar s))$ to the generalized graph $\Gamma_{\varphi_{12}\res[0,w_{12}]}$ of $\varphi_{12}\res [0,\walpbet]$ (for all $\bar s$ where it exists), see Figure \ref{fig:normaltriangle}. Without 
loss of generality we may assume $\Gamma_1 =\Gamma_{\varphi_{12}\res[0,w_{12}]}$.  
We start to show that $\phialpbet$ cannot have too negative slope, otherwise $\Gamma_1$ loses graphicality with respect to $\overline{\gam \alp}$.

\smallskip
\textit{Step 1.}
We claim that 
\begin{equation*}
\phialpbet^\prime\res [0,\walpbet] \geq\frac{\beta}{\alpha}
\end{equation*}
in the sense of measures, {\it i.e.}, 
\begin{equation}
\phialpbet^\prime(B)\geq\frac{\beta}{\alpha} \mathcal{L}^1(B), \qquad  \forall B\subseteq [0,\walpbet]{~\rm  Borel~ set. } \label{eqn:graphicalityphi}
\end{equation} 

From the graphicality with respect to  $\overline{\alp \wgamalp}$ we have, for all $\bar s$ where $ \nualp(\bar s)$ exists,
\begin{align}
\nualp(\bar s)\cdot\ngamalp\leq0. \label{5.41}
\end{align}
Set 
\begin{align*}
I^r:=\{\bar s\in [0,\walpbet]:\nualp(\bar s) {\rm ~is~defined, ~and ~}\nualp^2(\bar s)>0\},\\
I^s:=\{\bar s\in [0,\walpbet]:\nualp(\bar s) {\rm ~is~defined, ~and ~}\nualp^2(\bar s)=0\},
\end{align*}
note that $\nu(\bar s)=\frac{1}{\sqrt{1+(\dot{\varphi}_{12}(\bar s))^2}}(-\dot{\varphi}_{12}(\bar s),1)$ for any $\bar s\in I^r$. From \cite[Thm. 7 p. 301 and Thm. 5 p. 379]{Giq}, we have
\begin{equation}
\dot{\varphi}_{12}d \bar s={\varphi '}_{12} \res I^r, \qquad \dot{\varphi}_{2}^{(j)}+\dot{\varphi}_{12}^{(c)}={\varphi '}_{12} \res I^s.
\end{equation}

From \eqref{5.41} it follows that 
\begin{equation}
\nualp(\bar s)=(-1,0)\qquad \forall \bar s\in I^s
\qquad {\rm and} \qquad \dot{\varphi}_{12}(\bar s)\geq \frac{\beta}{\alpha} \qquad \forall \bar s\in I^r .
\label{eqn:phidotg}
\end{equation}
From \cite[Thm. 4, p. 378]{Giq} we have 
\begin{equation*}
\dot{\varphi}_{12}^{(j)}+\dot{\varphi}_{12}^{(c)}=-\nu^1\vert \mu \vert \res I^s=\vert \mu \vert \res I^s,\label{eqn:positivesinguler}
\end{equation*}
where $\mu:=(\phialpbet^\prime, -\mathcal{L}^1)=(-\nualp^1,-\nualp^2)|\mu|$ 
and the second equality follows from the first formula in \eqref{eqn:phidotg}.

For any Borel set $B\subseteq [0,\walpbet]$ we deduce
\begin{equation*}
\phialpbet^\prime(B)=\int_{B}\dot{\varphi}_{12}d\bar s+\dot{\varphi}_{12}^{(j)}(B)+\dot{\varphi}_{12}^{(c)}(B)\geq \frac{\beta}{\alpha}~\mathcal{L}^1(B)+\vert \mu \vert \res I^s(B)\geq \frac{\beta}{\alpha}~\mathcal{L}^1(B).
\end{equation*}

\smallskip

\textit{Step 2.}
Given $\epsilon \in (0,1)$, we choose $n=n(\epsilon)\in \N$ and points 
$$\xi_0=0< \xi_1<\cdots<\xi_{n-1}<\xi_n =\walpbet,$$
such that each $\xi_i,~i\in\{1,\cdots,n-1\}$, is a point of 
continuity of $\phialpbet$, and if we define $\fe \in \Lip([0,\walpbet])$ as the piecewise linear interpolation with 
$$\fe (\xi_i)=\phialpbet(\xi_i),\qquad i=0,\cdots,n, $$
then $$ \| \fe -\phialpbet\|_{L^1((0,\walpbet))}<\epsilon.$$

The graph of $\fe$ may still have vertical parts over $\overline{\wgamalp \alp}$. 
\begin{figure}[h!] 
  \begin{center}
   \includegraphics[scale=0.5]{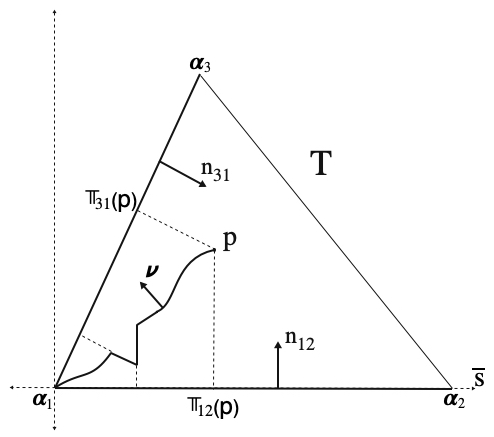}
           \caption{For convenience we choose $\alp=0$.} 
           \label{fig:normaltriangle}
  \end{center}
\end{figure}
Indeed from \cite[Theorem 3.30]{AmFuPa:00}, and the fact that $\xi_i$ are 
continuity points of $\phialpbet$ we have
\begin{equation}\label{eqn:fe}
({\fe})^\prime(\xi) =\frac{\phialpbet(\xi_i)-\phialpbet(\xi_{i-1})}{\xi_{i}-\xi_{i-1}}= \frac{\phialpbet^\prime((\xi_{i-1},\xi_{i}))}{\xi_{i}-\xi_{i-1}}\geq \frac{\beta}{\alpha},\qquad \xi \in(\xi_{i-1},\xi_{i}),
\end{equation}
and equality may hold, hence the graph of $\fe$ over $\overline{\pi_{31}(p)\alp}$ may have finitely many vertical parts. 
It is now sufficient to repeat the argument with $\fe$ in place of $\phialpbet$, choosing a suitable partition of $[\wgamalp,\ellgamalp],$ so to ensure that (out of finitely many points) $$ 
 ({\fe})^{\prime}> \frac{\beta}{\alpha}\qquad {\rm on~~}[0,\walpbet]. $$

In this way $\fe$ is a Lipschitz graph also with respect to $\overline{\pi_{31}(p)\alp}$.

\smallskip
\textit{Step 3.}
We have 
\begin{equation}
\begin{split}
&\Hone(\Gamma_{\fe})=\sum_{i=1}^n\int_{\xi_{i-1}}^{\xi_i}\sqrt{ 1+\vert {\fe} ^\prime (s) \vert^2 }ds=\sum_{i=1}^n\vert (\xi_i,{\fe} (\xi_i))-(\xi_{i-1}, {\fe}  (\xi_{i-1}))\vert \\ 
\leq &\sup \left\{  \sum_{i=1}^m\vert (\eta_i,{\phialpbet} (\eta_i))-(\eta_{i-1}, {\phialpbet}  (\eta_{i-1}))\vert: m\in \N, \eta_0=0< \eta_1<\cdots<\eta_{m-1}<\eta_m =\walpbet \right\}\\
=& \int_{[0,\walpbet]}\vert \Phi'_{12} \vert =\Hone({\Gamma_{\varphi_{12}\res [0,\walpbet]}}),
\end{split}\label{eqn:Honesup}
\end{equation}
where $ \Phi_{12}\in \BV([0,\walpbet];\R^2)$ is defined as $\Phi_{12}(\xi):=(\xi,\phialpbet(\xi))$, and the last equality follows from \cite[(3.24), p.136]{AmFuPa:00}.

 \textit{Step 4.}
  Define $$\Gammaalp^n:= {\Gamma}_{\fe}.$$ 
Similarly we define $\Gammabet^n$ and $\Gammagam^n$, and we set $\Gamma^n_{ij}:=\Gamma^n_{i}\cup\Gamma^n_{j}$. 
Then $\Gamma^n:=(\Gamma_{1}^n,\Gamma_{2}^n,\Gamma_{3}^n)$, 
satisfies the required properties.   
\end{proof}

\begin{Proposition}[\textbf{Uniform estimate of the length}]
\label{pro:uniform_estimate_on_the_length}
There exists $c>0$ such that for all $\elementXclip \in X$ we have 
\begin{equation}
\Hone(\Gammaalpij)\leq c,\qquad  ij\in \{12,23,31\}.\label{eqn:lessc}
\end{equation}
\end{Proposition}
\begin{proof}
Let $\elementXclip \in \Xc$ be a connection through $\genericp\in \targettriangle$. Without 
loss of generality we may assume that $\genericp\not= \alp$. From \eqref{eqn:Honesup} we have 
\begin{equation}
\Hone(\Gammaalp)=\sup \left\{  \sum_{i=1}^m\vert (\eta_i,{\phialpbet} (\eta_i))-(\eta_{i-1}, {\phialpbet}  (\eta_{i-1}))\vert:m\in \N, \eta_0=0< \eta_1<\cdots<\eta_{m-1}<\eta_m =\walpbet \right\}.\label{eqn:gammalpsup}
\end{equation}
Choose a partition
$$ \xi_0=0< \xi_1<\cdots<\xi_{\hkindex-1}<\xi_\hkindex=\walpbet.$$
Let $\Gamma_1^\hkindex$ be the piecewise linear interpolation connecting $(\xi_{i-1},\phialpbet(\xi_{i-1}))$ and $(\xi_i,\phialpbet(\xi_{i})),$ $i\in\{1,\cdots,\hkindex\}$. The unit tangent to $\Gamma_1^\hkindex$ is enclosed in the angle formed by $n_{12}$ and $n_{31}$, the unit normals to $\overline{\alp\bet}$ and $\overline{\gam \alp}$ (due to the graphicality condition with respect to $\overline{\alp\bet}$ and $\overline{\gam \alp}$), see Figure \ref{fig:uniformlengthtriangletwo}. It follows that $\Gamma_1^\hkindex$ is the graph of a function $\phi_{12}^\hkindex$ over the segment $\overline{\alpha_1p}$. Fix a Cartesian coordinate system in which the 
$t$-axis is the line $\overline{\alp \genericp}$ and the origin is $\alpha_1$. 
 \begin{figure}
\centering
\begin{subfigure}{.5\textwidth }
  \centering
  \begin{center}
     \includegraphics[scale=0.45]{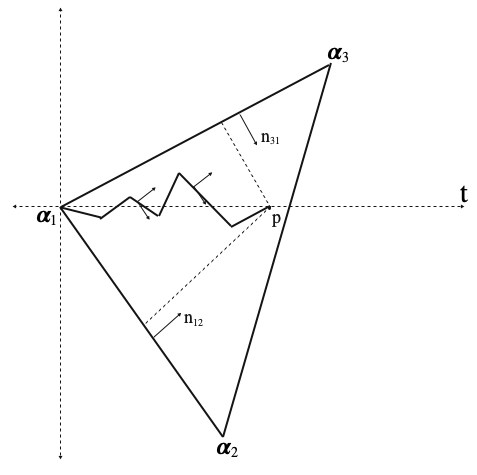}
   \end{center}
 \caption{$\targettriangle$ with angles less than or equal to $\frac{\pi}{2}$. }     \label{fig:uniformlengthtriangle}
  \end{subfigure}%
\begin{subfigure}{.5\textwidth} 
  \centering
  \begin{center}
       \includegraphics[scale=0.4]{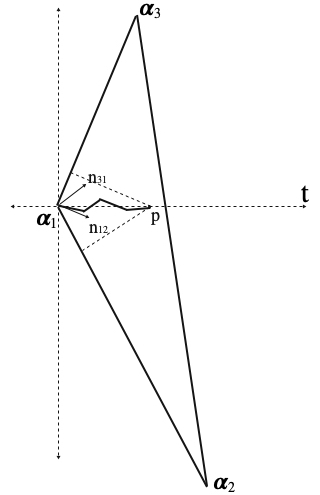}
  \end{center}
  \caption{$\targettriangle$ with an angle greater than $\frac{\pi}{2}$.}
  \label{fig:uniformlengthtrianglecopy}
\end{subfigure}
\caption{}\label{fig:uniformlengthtriangletwo}
\end{figure}
For any $t\in [0, \vert \alp -\genericp \vert]$ 
(up to a finite set) let $\tau(t)$ be the unit tangent to $\Gamma_1^\hkindex$ at $(t, \phi^\hkindex_{12}(t))$ and let $n=(1,0)$ and $n^\perp=(0,1)$.
Hence ${\phi^\hkindex_{12}}'=\frac{\tau\cdot n^\perp}{\tau\cdot n}$ satisfies
\begin{equation*}
c_1^-:=\frac{\ngamalp\cdot n^\perp}{\ngamalp\cdot n}\leq {\phi^\hkindex_{12}}' \leq \frac{\nalpbet\cdot n^\perp}{\nalpbet\cdot n}=:c_1^+.
\end{equation*}
{Note that one between $|c_1^-|$ and $|c_1^+|$ might be $+\infty$, since one of the sides $\overline{\alp\bet}$ or $\overline{\gam \alp}$ can be horizontal (this happens only if the point $p$ is on one side of the triangle). However 
we always have that $c_1^-\leq0$, $c_1^+\geq0$.
Furthermore, when the angle $\hat\alpha_1$ in $\alpha_1$ is less or equal to $\frac{\pi}{2}$, it follows that  $\tilde c_1:=\min \{|c_1^-|,|c_1^+|\}\leq |\tan(\frac{\pi}{2}-\frac{\hat\alpha_1}{2})|$. In the case that $\hat\alpha_1>\frac{\pi}{2}$, thanks to the fact that $p\in T_{int}$, we have $\max \{|c_1^-|,|c_1^+|\}\leq|\tan(\pi-\hat\alpha_1)|$. Thus the only difficulty to prove that the length of $\Gamma_1^\hkindex$ is controlled when $\hat\alpha_1\leq\frac{\pi}{2}$. So let us assume this and in addition that $|c^-_1|=\tilde c_1$ (the other case is similar).
Since $\phi^\hkindex_{12}(\vert \alp -\genericp \vert)=\phi^\hkindex_{12}(0)=0$ we have 
\begin{align*}
0={\phi^\hkindex_{12}}'([0,\vert \alp -\genericp \vert]) =({\phi^\hkindex_{12}}')^+([0,\vert \alp -\genericp \vert])-({\phi^\hkindex_{12}}')^-([0,\vert \alp -\genericp \vert]),
\end{align*}
where $({\phi^\hkindex_{12}}')^+$ and $({\phi^\hkindex_{12}}')^-$ are the positive and negative parts of the measure ${\phi^\hkindex_{12}}'=\dot\phi^\hkindex_{12} dt$,
thus we estimate
\begin{align}\label{5.15}
|\mathcal H^1(\Gamma_1^\hkindex)|&=
\int_0^{\vert \alp -\genericp \vert}\sqrt{1+\dot\phi^\hkindex_{12}(t)^2}dt\nonumber\\
&\leq\vert \alp -\genericp \vert+|{\phi^\hkindex_{12}}'|([0,\vert \alp -\genericp \vert])=\vert \alp -\genericp \vert+2({\phi^\hkindex_{12}}')^-([0,\vert \alp -\genericp \vert])\nonumber\\
&\leq \vert \alp -\genericp \vert+2\vert \alp -\genericp \vert\tilde c_1.
\end{align}
Defining $c_1$ as the right-hand side of the last inequality we see that $c_1$ is a positive constant depending only on the geometry of $\targettriangle$.

 From the \eqref{5.15} and \eqref{eqn:gammalpsup} it then follows
\begin{equation}
\Hone(\Gammaalp)\leq c_1.
\end{equation}
Similarly we may show that $\Hone({\Gammabet})\leq c_2$ and $\Hone({\Gammagam})\leq c_3$ for $c_2, c_3>0$ depending only on $\targettriangle$. This proves \eqref{eqn:lessc} with $c=c_1+c_2+c_3$.
}
\end{proof}

The next lemma shows 
continuity of the sum of the three areas of area minimizing surfaces defining $\AG$ in \eqref{eqn:AG}, with respect to the $\Lone$ convergence of the traces in $\targettriangle$. 
\begin{Proposition}[\textbf{Continuity of $\AG$}]\label{Lem:AGcontinuous}
Let $\elementXclip \in X$, and let $\{\elementXclip^n\}\subset X$ be a sequence converging to $\elementXclip$ in $X$.
Then
 \begin{equation}
\lim_{n\to +\infty}\AG (\elementXclip^n) = \AG(\elementXclip).\label{eqn:4.61}
\end{equation}
 \end{Proposition}
\begin{proof}
Since $\elementXclip \in X$ and $\{\elementXclip^n\}\subset X$ we have $\phiij \in {~\rm BV}([0,\ellalpij])$ and $ \{\phiij^n\} \subset {~\rm BV}([0,\ellalpij]) $ where  $\phiij:=\phiij(\pGammaalpij)$, $\phiij^n:=\phiij(\pnGammaalpij)$.

Hence from \eqref{eq:extension_of_varphi} and Section \ref{tenda} it follows that there exist $\widehat \tenda_ {{ij}},~\widehat \tenda^n_{{ij}} \in W^{1,1}(\widehat\rectangle_{{ij}})$  such that
 \begin{align}
\begin{split}
2{\mathfrak A}^n_{i j}=&\int_{\widehat\rettangolo_{i j}} \sqrt{1 + \vert\grad \widehat \tenda^n_ {{ij}}\vert^2}~dsdt\\
=&\min \big\{ \int_{\widehat\rettangolo_{i j}} \sqrt{1 + \vert D f\vert^2}~+\int_{\partial \widehat \rectangle_{{ij}}}|f-\phiij^n|d\H^{1}: f \in {\rm BV}(\Ball),~f=\phiij^n {~\rm on~}\Ball\setminus \widehat \rectangle_{{ij}}\big\},\end{split} \label{minphin}\\
\begin{split}
2{\mathfrak A}_{i j}=&\int_{\widehat\rettangolo_{i j}} \sqrt{1 + \vert\grad \widehat \tenda_ {{ij}}\vert^2}~dsdt \\
=&\min \big\{ \int_{\widehat\rettangolo_{i j}} \sqrt{1 + \vert Df\vert^2}~+\int_{\partial \widehat \rectangle_{{ij}}}|f-\phiij| d\H^{1}: f \in {\rm BV}(\Ball),~f=\phiij {~\rm on~}\Ball\setminus \widehat \rectangle_{{ij}}\big\},\end{split}\label{minphi}
\end{align}
where we recall that $\widehat\rectangle_{{ij}}$ is the double rectangle defined in \eqref{eqn:Rhat} and $\phiij,~\phiij^n$ are extended on a disk $\Ball$ containing $\widehat\rectangle_{{ij}}$ as in Section \ref{tenda}. 

 Define $\widetilde \tenda^n_{{ij}}$ and $\widetilde \tenda_{{ij}}$  as 
  \begin{equation*}
  \widetilde \tenda_{{ij}}^n := 
  \begin{cases*}
  \widehat \tenda^n_{{ij}} {\rm~in~}\widehat \rectangle_{{ij}},\\
  \phiij  {\rm~in~}\Ball\setminus \widehat \rectangle_{{ij}},
  \end{cases*} 
 \qquad \widetilde \tenda_{{ij}}:=
  \begin{cases*}
  \widehat \tenda_{{ij}} {\rm~in~}\widehat \rectangle_{{ij}},\\
  \phiij^n  {\rm~in~}\Ball\setminus \widehat \rectangle_{{ij}},
  \end{cases*}
 \end{equation*}
  so that $\widetilde \tenda^n_{{ij}},~\widetilde \tenda_{{ij}} \in {\rm BV}(\Ball)$.  
 Since $\widetilde \tenda^n_{{ij}}$ is competitor in \eqref{minphi} and $\widetilde \tenda_{{ij}}$ is competitor in \eqref{minphin}
 we have, recalling also the discussion leading to \eqref{widehat_tenda}, 
  \begin{eqnarray*}
2{\mathfrak A}_{i j}\leq \int_{\widehat\rettangolo_{i j}} \sqrt{1 + \vert\grad \widehat \tenda^n_{{ij}}\vert^2}~dsdt+\int_{\partial \widehat \rectangle_{ij}}|\phiij^n-\phiij| d\H^{1}= 2{\mathfrak A}^n_{i j}+\int_{\partial \widehat \rectangle_{ij}}|\phiij^n-\phiij| d\H^{1},\\
2{\mathfrak A}^n_{i j}\leq \int_{\widehat\rettangolo_{i j}} \sqrt{1 + \vert\grad \widehat \tenda_{{ij}}\vert^2}~dsdt+\int_{\partial \widehat \rectangle_{ij}}|\phiij-\phiij^n| d\H^{1}= 2{\mathfrak A}_{i j}+\int_{\partial \widehat \rectangle_{ij}}|\phiij-\phiij^n| d\H^{1}. 
\end{eqnarray*}
 Thus 
 \begin{equation}
 \vert 2{\mathfrak A}^n_{i j}-2 {\mathfrak A}_{i j}\vert  \leq \int_{\partial \widehat \rectangle}\vert\phiij^n-\phiij\vert d\H^{1}.\label{eqn:AGleq}
 \end{equation}
Recall that $\tenda^n_{{ij}}$ (resp. $\tenda_{{ij}}$) is the restriction of $\widehat \tenda^n_{{ij}}$ (resp. $\widehat \tenda_{{ij}}$) to $\rectangle_{{ij}}$. Hence, from \eqref{eqn:AG}, \eqref{Def:XConv} and \eqref{eqn:AGleq},
\eqref{eqn:4.61} follows.
\end{proof}

\begin{Corollary} \label{cor:density}
We have
\begin{equation}\label{eqn:infXcequalinfXclip}
\inf \left\{{\AG}(\elementXclip): \elementXclip \in \Xc\right\}=\inf \left\{{\AG}(\elementXclip): \elementXclip \in \Xclip\right\}.
\end{equation}
\end{Corollary}


\subsection{Compactness of the class $\Xc$}
The aim of this section is to show that the infimum in \eqref{eqn:infXcequalinfXclip} is attained. To do this we need the following result.
\begin{Theorem}[\textbf{Compactness}]\label{teo:closureXclip}
 Any sequence $\{\elementXclip^n\}\subset \Xc$ admits a subsequence converging in $X$ to some $\elementXclip \in \Xc$. 
\end{Theorem}

\begin{Remark}\rm
 In Definition \ref{ConvinX} 
it is required convergence of $\{\Gamma^n\}$ to $\Gamma$ in $L^1$. For this 
reason, if $\Gamma^n$ has target triple point $p_n$,  it is not guaranteed 
that the point $b:=\lim_{n\rightarrow +\infty}p_n$ (it exists up to subsequences)
 still belongs to $\Gamma_{ij}$ for all ${ij}$,  see Figures \ref{fig2comp} 
and \ref{fig:caseiva}. As a consequence, if $\{\elementXclip^n\}$ converges to $\elementXclip$ it is not true, in general, that $p_n\rightarrow p$, where $p$ is the target triple point of $\elementXclip$. 
\end{Remark}

\begin{proof}
Let $\{\elementXclip^n\}\subset \Xc$ and $\phialpij^n=\phialpij(\pnGammaalpij),~ {ij} \in \{12,23,31\}$. 
From Proposition \ref{pro:uniform_estimate_on_the_length} 
$\{\phialpij^n\}$ is 
uniformly bounded in $\BV([0,\ellalpij])$ for any ${ij} \in \{12,23,31\}$. Thus, up to a not relabelled subsequence, there exists $\phialpij \in \BV([0,\ellalpij])$ such that
\begin{align}
 &\phialpij^n \to \phialpij \qquad{\rm ~ in ~ } \Lone((0,\ellalpij)) {\rm~ and ~pointwise~ \rm~a.e.},\label{eqn: strong_convergence_L1}\\
 &(\phialpij^n)^\prime \rightharpoonup \dirdatum_{ij}^\prime \qquad {\rm weakly}^{*} {\rm~as~measures}\label{eqn: convergence_measure}.
 \end{align} 
We shall adopt our usual convention \begin{equation}
\phialpij^n(0_-)=\phialpij(0_-)=\phialpij^n(\ellalpij_+)=\phialpij(\ellalpij_+)=0,
\;\;\varphi_{ij}^n={\varphi_{ij}^n}_+,
\;\;\varphi_{ij}={\varphi_{ij}}_+,
\;\;\varphi=\varphi_+.
\label{eqn:phiij0}
\end{equation}
Denote by $\Gammaalpij \subset \R^2$ the limit graph over (the closed segment) $\overline{\alpi \alpj}$ that we identify with the generalized graph of $\phialpij$ over $[0,\ellalpij]$. 
Since $\targettriangle$ is closed and convex we have 
 $\Gamma_{ij}\subset T$; moreover, by construction, $\alpha_i$ and $\alpha_j$ are the endpoints of $\Gamma_{ij}$. 
Notice that if we assume that $\targettriangle$ is acute,
this excludes the presence of 
vertical parts over its vertices. 

It remains to prove
that the three obtained curves $\Gamma_{ij}$, ${ij} \in \{12,23,31\}$, form
a BV connection; in particular that they intersect mutually in a unique
well-defined point.

\smallskip
We claim that 
\begin{itemize}
 \item[] there exists a unique $p \in \bigcap_{ij}\Gammaalpij$ that divides each $\Gammaalpij$ into two curves $\Gamma_{ij}^l$ and $\Gamma_{ij}^r$ 
such that $$\Gamma_{ij}^l=\Gamma_{ki}^r,\qquad   {ij}, \ki \in \{12,23,31\},~{ij} \not = \ki.$$
\end{itemize}

Let us denote by ${\tilde \varphi_{ij}^n}$ 
the extension to $\R$ of the function $\phialpij^n$ 
vanishing
in $(-\infty,0)\cup (\ellalpij,+\infty)$. 
Similarly $\tilde \varphi_{ij}$ is the extension of $\varphi_{ij}$
vanishing
in $(-\infty,0)\cup (\ellalpij,+\infty)$.
Consider the sequence $\{[\![S\mathcal G_{\tilde \varphi_{ij}^n}]\!]\}_n\subset \mathcal D_2(\R^2)$ of 2-currents regarded in $\R^2$ and the 2-current $[\![ S\mathcal G_{\tilde \varphi_{ij}}]\!]$. Their boundaries are the currents carried by the graphs of  $\tilde \varphi_{ij}^n$ and  $\tilde\varphi_{ij}$, as defined in Theorem \ref{thm: boundary currents}. 
The 1-currents carried by the graph of $\varphi_{ij}^n$ and $\varphi_{ij}$, by convention \eqref{eqn:phiij0}, coincide with the restrictions of $\partial[\![S\mathcal G_{\tilde \varphi_{ij}^n}]\!]$ and $\partial[\![ SG_{\tilde \varphi_{ij}}]\!]$ to the closed set $[0,\ellalpij]\times \R$. 
Namely, if we denote by
\begin{align*}
 \currentijn:=\partial[\![S\mathcal G_{\tilde \varphi_{ij}^n}]\!]\res [0,\ellalpij]\times 
\R,\;\;\currentij:=\partial[\![ SG_{\tilde \varphi_{ij}}]\!]\res [0,\ellalpij]\times \R,
\end{align*}
then 
\begin{align}\label{1000}
 \currentijn=\partial[\![S\mathcal G_{\tilde \varphi_{ij}^n}]\!]-\mathcal L_{ij}\text{ and }\currentij=\partial[\![ SG_{\tilde \varphi_{ij}}]\!]-\mathcal L_{ij},
\end{align}
where $\mathcal L_{ij}$ is the $1$-current given by integration over the two 
halflines $(-\infty,0)\times\{0\}\cup(\ellalpij,+\infty)\times\{0\}$. 
The curves  $\Gamma_{ij}^n$ and $\Gamma_{ij}$ coincide with the support of $\currentijn$ and $\currentij$, respectively.

We now prove our claim in three steps.

\smallskip
{\it Step 1.}  The currents $\currentijn$ converge 
(up to a not relabelled subsequence) weakly in the sense of currents to {$\currentij$}, {\it i.e.},
\begin{equation}
\currentijn(\omega) \to \currentij (\omega) \qquad \forall  \omega \in \DD^1(\R^2). \label{eqn:convtoGammaij} 
\end{equation}
Moreover
\begin{equation}
\Hone(\Gammaalpij)\leq c,\label{eqn:Honec} 
\end{equation}
where $c>0$ is the constant in \eqref{eqn:lessc}.

Indeed, 
thanks to \eqref{eqn: strong_convergence_L1}, 
the characteristic functions $\chi_{S\G_{\phiij^n}}$ converge to $\chi_{S\G_{\phiij}}$  in $L_{{\rm loc}}^1(\R^2)$, hence
\begin{align*}
  [\![S\mathcal G_{\tilde \varphi_{ij}^n}]\!]\rightharpoonup  [\![S\mathcal G_{\tilde \varphi_{ij}}]\!]\text{ weakly as currents,}
\end{align*}
 since
\begin{align*}
 \int_{S\G_{\tilde \varphi_{ij}^n}}\hat\omega(s,t)dsdt\rightarrow\int_{S\G_{\tilde\varphi_{ij}}}\hat\omega(s,t)dsdt
\qquad
\forall
\hat\omega\in\mathcal D^2(\R^2).
\end{align*}
This implies 
$$
\partial[\![S\mathcal G_{\tilde \varphi_{ij}^n}]\!]\rightharpoonup  \partial[\![S\mathcal G_{\tilde \varphi_{ij}}]\!]\text{ weakly in the sense of currents,}
$$
and \eqref{eqn:convtoGammaij} follows from \eqref{1000}. 

Finally \eqref{eqn:Honec} follows from Lemma \ref{pro:uniform_estimate_on_the_length} and the 
weak lower semicontinuity of the mass of currents, and the proof of 
step 1 is concluded.

\smallskip
It is not restrictive to assume that
$\wij^n=|\alpi-\pi_{{ij}}(\genericpn)|$ is a point of continuity of $\phiij^n$
 for all $n \in \N$ and all ${ij}\in \{12,23,31\}$. 
Indeed given a sequence $\{\elementXclip^n\}\subset \Xc$ converging to $\Gamma$, 
from Lemma \ref{Lem:density}
for all $n$ 
we can  
assign a sequence $\{\elementXclip^{m,n}\}\subset \Xclip$ such that 
$\elementXclip^{m,n}\rightarrow\Gamma^n$ as $m\rightarrow+\infty$. Thus by a diagonal argument, we find a sequence $\{\elementXclip^{m(n),n}\}\subset \Xclip$ which tends to $\Gamma$ and satisfies 
the above requirement
(we can also assume that 
$\Gamma^n$ is Lipschitz, but this will not
be needed in the proof).

Without loss of generality 
(up to a not relabeled subsequence) we may 
further assume
$$
\genericpn \to b\in \targettriangle,
$$
$\{\wij^n\}$ is a monotone sequence, and
$$
\wij^n \to \wij:=\vert \alpha_i - 
\pi_{ij}(b)\vert, \qquad ij\in \{12,23,31\}.
$$
Before passing to the second
step, it is convenient to divide the target triangle $\targettriangle$ into 
various regions.

Assume first that $\targettriangle$ is acute.
The point $b$, together with the heights 
\begin{equation}
\label{eqn:Lij}
\Lij:=\overline{\pi_{ij}( b)b},\qquad {ij} \in \{12,23,31\},
\end{equation} 
divides $\targettriangle$ 
into three regions $\Ri,~i\in\{1,2,3\}$, 
as shown in Figure $\ref{fig1comp}$; precisely,
if $\overline{\mathcal P}_i$ denotes the closed region  enclosed by $h_{ij}$, $h_{ki}$, $\overline{\alpi \pi_{ij}{(b)}}$ 
and $\overline{\alpi \pi_{ki}{(b)}}$, then ${\mathcal P_i}$ is defined by 
\begin{equation}
\mathcal P_i:= \overline{\mathcal P}_i\setminus (h_{ij}\cup h_{ki}),\qquad i=1,2,3.\label{eqn:Rij}
\end{equation}
Similarly we define $\Lij^n$ and $\Ri^n$ by replacing $b$ with $\genericpn$ in \eqref{eqn:Lij} and \eqref{eqn:Rij}. 

Assume now that $\targettriangle$ is not acute. Without 
loss of generality we may assume that the angle at $  \alp$ is greater than $\frac{\pi}{2}$. The only difference here is with the definition of $\Ralp^n$ and $\Ralp$ since each $\Gammaalpij$ has to satisfies the graphicality condition with respect to $\overline{\alpha_i \alpha_j}$; hence we define $\overline{\Ralp}$ as the closed quadrilateral bounded by  $\Lalpbet$, $\Lgamalp$, $m_{12}$ and $m_{31}$, where $m_{12}$ and $m_{31}$ are the normals to $\overline{\alp \bet}$ and $\overline{\gam \alp}$, respectively, passing through $\alp$ (see Figure \ref{fig:T180partitions}).
Similarly we define $\overline{\Ralp}^n$. Finally we set $\Ralp:=\overline{\Ralp}\setminus(\Lalpbet\cup\Lgamalp)$ and $\Ralp^n:=\overline{\Ralp}^n\setminus(\Lalpbet\cup\Lgamalp)$.

\smallskip
{\it Step 2.} We will prove that we can decompose $\Gamma_{12}\cup \Gamma_{23}\cup \Gamma_{31}$ as three currents meeting at a point $b$.

It is 
easy to see that the sets $\overline{\mathcal P}_i^n$ are converging to $\overline{\mathcal P}_i$ with respect to the Hausdorff distance. It is not true in general that $\Gamma_{ij}^n$ is converging to $\Gamma_{ij}$ with respect to the Hausdorff distance (see Figure \ref{fig2comp}); however, since
\begin{align}
\bGammaalpi=\pnGammaalpij\cap\overline{\Ri^n}=\pnGammaalpki\cap\overline{\Ri^n},\qquad
\bGammaalpi \subset \overline{\mathcal P}_i^n,\label{eqn:star}
\end{align}
for all $ij,ki\in\{12,23,31\}, ij \not = ki$, it is readily seen that 
\begin{align}
&\Gamma_{ij}\subset\targettriangle\setminus\Rk,\nonumber\\
&\currentij=\currentij\res\Ri+\currentij\res\Rj+\currentij\res\Lki+\currentij\res\Lij+\currentij\res\Ljk.\label{eqn:Gammaijdevision}
\end{align} 

\begin{figure}
\centering
\begin{subfigure}{.53\textwidth }
  \centering
  \begin{center}
     \includegraphics[scale=0.4]{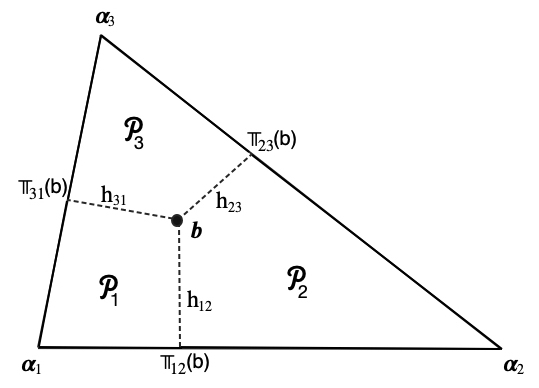}
   \end{center}
 \caption{Partitions of $\targettriangle$ into $\mathcal P_1$, $\mathcal P_2$, $\mathcal P_3$ and three segments.} 
           \label{fig1comp}
  \end{subfigure}%
\begin{subfigure}{.47\textwidth} 
  \centering
  \begin{center}
       \includegraphics[scale=0.36]{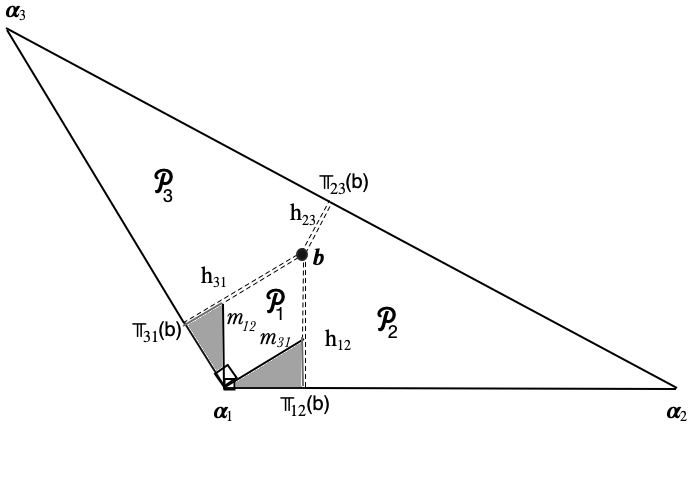}
  \end{center}
  \caption{ Partitions of $\targettriangle$ minus the two grey triangles intointo $\mathcal P_1$, $\mathcal P_2$, $\mathcal P_3$ and three segments for a non acute triangle.}\label{fig:T180partitions}
\end{subfigure}
\caption{}
\end{figure}

For $i\in\{1,2,3\}$, the 
integral 1-current $\currentib=\currentijn\res[0,w_{ij}]\times\R$ has boundary  $\delta_{\genericpn}-\delta_{\alpi}$ in $\R^2$.  By the 
compactness theorem for integral currents \cite[Theorem 2, p.141]{Giq} there exists an integral current $\Tli \in \mathcal D_1(\R^2)$, $i=1,2,3$,  such that, up to a not relabeled subsequence,
\begin{equation}
\currentib(\omega)\to\Tli({\omega}) \qquad \forall \omega\in \DD^1(\R^2).\label{eqn:convtoTi}
\end{equation}
Clearly
 \begin{equation}
 \partial\Tli =\delta_b-\delta_{\alpi}.\label{eqn:boundaryTi}
 \end{equation}
From \eqref{eqn:star} and thanks to the convergence of $\overline{\mathcal P_i^n}$ to $\overline{\mathcal P_i}$ with respect to the Hausdorff distance, we infer  
\begin{equation}
\spt\Tli\subset \overline\Ri,\qquad{\rm hence~} \Tli=\Tli\res\Ri+\Tli\res\Lij+\Tli\res\Lki,\label{eqn:Tidevision}
\end{equation} 
where ${ij},\ki\in \{12,23,31\}$.
Note that $\Tli$ is not necessarily equal to $\currentij \res \overline \Ri$, due to a possible cancellation of a vertical part over $\pi_{ij}(b),~ij\in\{12,23,31\}$ (that is, on $h_{ij}$), see Figure \ref{fig2comp}. 
However from $\currentijn=\currentib+[\![\Gamma_j^n]\!]$ and \eqref{eqn:convtoTi} we have 
\begin{equation}
\currentij =\Tli-\Tlj, \qquad {ij}\in \{12,23,31\},\label{eqn:TiminusTj}
\end{equation}
as currents in $\R^2$. Notice that $\Tli$ and $\Tlj$ have multiplicity one, and in \eqref{eqn:TiminusTj} they contribute with opposite orientation. This allows, if necessary, to identify $\Tli,~i=1,2,3,$ with its support.  
\begin{figure}[h!] 
  \begin{center}
   \includegraphics[scale=0.5]{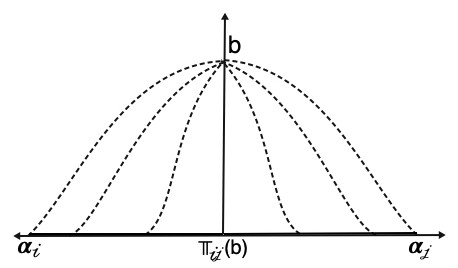}
           \caption{In dots a sequence of graphs $\currentijb$ of functions that pass through a fixed point $b\in T$. In bold the graph of the limit function (the horizontal segment) $\currentij$. The limit in the sense of currents of the left branches of the sequence $\{\Gamma_i^n\}$ is $\overline{\alpha_i \pi_{ij}(b)} \cup \overline{\pi_{ij}(b) b}$ while the limit of the right branches $\{\Gamma_j^n\}$ is $\overline{\pi_{ij}(b) b} \cup \overline{\pi_{ij}(b)\alpha_j }$ . } 
           \label{fig2comp}
  \end{center}
\end{figure}
Note also that $\Tli$ may have vertical part over $\alpi$, see Figure \ref{fig:T180}.
\begin{figure}[h!] 
  \begin{center}
   \includegraphics[scale=0.4]{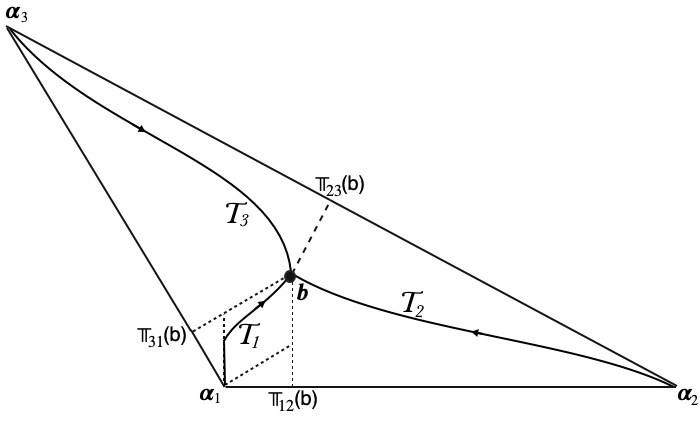}
    \caption{}
   \label{fig:T180}
  \end{center}
\end{figure}
Now, since $\currentib$ is Cartesian with respect to both the edges $\overline{\alpha_{i}\alpha_j}$ and $\overline{\alpha_{k}\alpha_i}$, from  \eqref{eqn:convtoTi} it follows that  $\Tli\res\Ri$ is part of two generalized graphs over the same edges, {\it i.e.}, 
\begin{equation}
\Tli\res \Ri=\currentij\res\Ri=-\currentki\res\Ri.\label{eqn:TiRi}
\end{equation}
Moreover, we infer that $\mathcal T_i$ cannot have vertical part over $h_{ij}$ and $h_{ki}$ at the same time; in other words once the current $\Tli$ touches one of the heights $\Lij$ or $\Lki$ it stays there until it reaches $b$, and $\Tli$ cannot have a nonempty support in more than one height, see Figures \ref{fig:caseii}-\ref{fig:caseiv}.
We conclude the following statement:

\begin{itemize}
 \item[(A)] The supports of the three currents $\mathcal T_i$, $i=1,2,3$, have as common point $b$. Moreover, if there are $i\neq j$ such that the supports of $\mathcal T_i$ and $\mathcal T_j$ intersect in a point different from $b$, then this intersection occurs on the mutual height $h_{ij}$. Finally, if the supports of $\mathcal T_i$ and $\mathcal T_j$ intersect on $h_{ij}$ outside $b$, then they intersect on a closed segment and the intersection of the supports of $\mathcal T_i$ and $\mathcal T_j$ with $\mathcal T_k$ is only the point $b$.   
\end{itemize}

\smallskip
\textit{Step 3.}
To conclude the proof of our claim 
we now analyse the possible cases arising from (A).

Case (i). Assume that the three supports of the currents $\mathcal T_i$, $i=1,2,3$, intersect only at the point $b$. This includes the case   $$\Tli \res \Lij = \Tlj \res \Lij=0{,\qquad\rm for~ all ~}{ij} \in \{12,23,31\},$$ as in Figure \ref{fig:casei}. But it may also happen that $\mathcal T_i$ has vertical part over $h_{ij}$, provided that $\mathcal T_j$ does not have vertical part over the same height (see for instance Figure \ref{fig:caseii}). 
In any case  we may set 
\begin{equation*}
p:=b,\qquad \Gamma_{ij}^l=\Tli,\qquad\Gamma_{ij}^r:=-\Tlj, \qquad{ij} \in \{12,23,31\},\label{eqn:casei}
\end{equation*}
where we have identified the currents $\mathcal T_i$ with their supports.
By \eqref{eqn:boundaryTi} and \eqref{eqn:TiRi}, the claim is achieved.

\begin{figure}
\centering
\begin{subfigure}{.5\textwidth }
  \centering
  \begin{center}
     \includegraphics[scale=0.45]{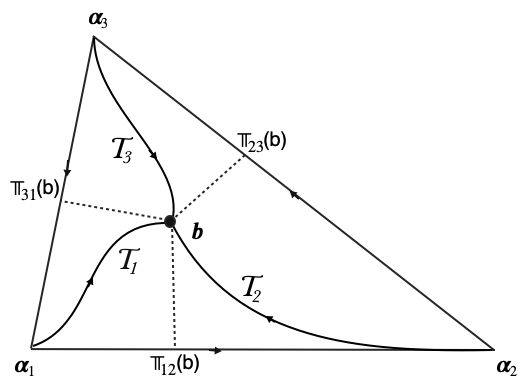}
   \end{center}
 \caption{} 
           \label{fig:casei}
  \end{subfigure}%
\begin{subfigure}{.5\textwidth} 
  \centering
  \begin{center}
       \includegraphics[scale=0.45]{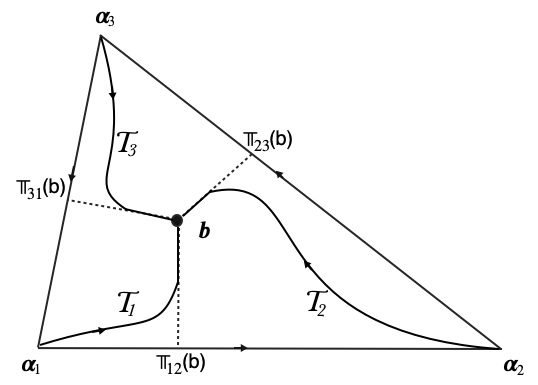}
  \end{center}
  \caption{}\label{fig:caseii}
\end{subfigure}
\caption{Case (i) of step 3 in the
 proof of Theorem \ref{teo:closureXclip}. }
\end{figure}

Case (ii). The second case which must be discussed is the one considering possible overlapping of the support of the currents $\mathcal T_{i}$. By condition (A)
such overlapping, giving rise to cancellations, can occur only on one height $h_{ij}$. Hence, assume  there exists one (and only one) ${ij} \in \{12,23,31\}$ such that $$\Tli \res \Lij \not= 0{\quad \rm and\quad }\Tlj\res\Lij \not=0.$$ Thus we have $\Tli\res\Lki=0$ and $\Tlj\res\Ljk=0$.

First assume that $\currentij\res\Lij=0$, {\it i.e.}, $\phiij$ is continuous at $\wij$. Then $\Tli\res\Lij=\Tlj\res\Lij$, see Figure \ref{fig:caseiva}. We set, identifying $\mathcal T_i$ with its support,
\begin{align*}
&p:=\phiij(\wij),&&&\\
&\Gamma_{ij}^l:=\Tli\res\Ri,&&\Gamma_{ij}^r:=\Tlj\res\Rj,&\\
&\Gamma_{jk}^l:=\Tlj\res\Ri,&&\Gamma_{jk}^r:=\Tlk\cup\Tlj\res\Lij,&\\
&\Gamma_{ki}^l:=\Tlk \cup \Tli\res\Lij,&&\Gamma_{ki}^r:=\Tli\res\Ri.&
\end{align*}
\begin{figure}
\centering
\begin{subfigure}{.5\textwidth }
  \centering
  \begin{center}
     \includegraphics[scale=0.45]{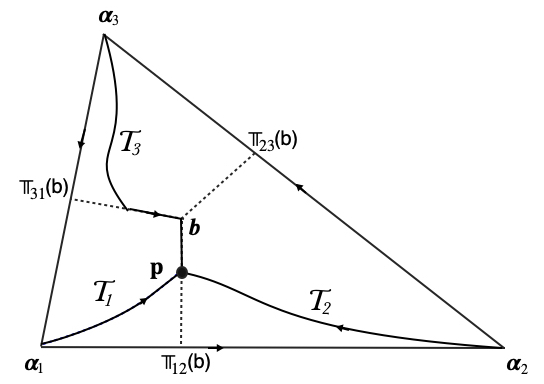}
   \end{center}
 \caption{$\T_1$ and $\T_2$ coincide only on $\overline{bp}$.} 
           \label{fig:caseiva}
  \end{subfigure}%
\begin{subfigure}{.5\textwidth} 
  \centering
  \begin{center}
       \includegraphics[scale=0.45]{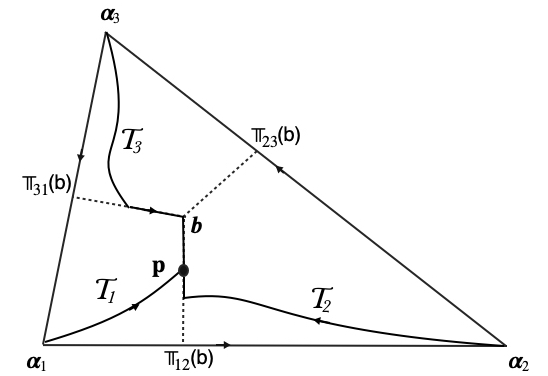}
  \end{center}
  \caption{$\T_1$ and $\T_2$ overlap on $\overline{bp}$ and $\T_1 \res \overline{bp} \subset \T_2 \res \overline{bp}$.}
  \label{fig:caseiv}
\end{subfigure}
\caption{Case (ii) of step 3 in the
 proof of Theorem \ref{teo:closureXclip}. }
\end{figure}
One checks  that the connection built above is a BV graph type connection, addressing the claim.

Now assume that $$\currentij\res\Lij\not=0,$$ {\it i.e.}, $\phiij$ jumps at $\wij$. Thus either $\spt\currentij\res\Lij \subseteq\spt\Tli\res\Lij$ or $\spt\currentij\res\Lij \subseteq\spt\Tlj\res\Lij$. Without 
loss of generality we may assume that $\spt \currentij\res\Lij \subseteq\spt\Tli \res\Lij$ hence 
$\Tli\res(\Lij\setminus \spt \currentij)=-\Tlj \res \Lij \not=0$  (note that $\spt\currentij\res\Lij=\{t\phiij(\wij{}_{+})+(1-t)\phiij(\wij{}_{-}):t\in[0,1]\}$).
We set
\begin{align*}
&p:=\phiij(\wij)=\varphi_{ij+}(\wij),&&&\\
&\Gamma_{ij}^l:=\Tli\res\Ri\cup(\Lij\cap \spt \currentij),&&\Gamma_{ij}^r:=\Tlj\res\Rj,&\\
&\Gamma_{jk}^l:=\Tlj\res\Ri,&&\Gamma_{jk}^r:=\Tlk\cup\Tlj\res\Lij,&\\
&\Gamma_{ki}^l:=\Tlk \cup \Tlj\res\Lij,&&\Gamma_{ki}^r:=\Tli\res\Ri\cup(\Lij\cap \spt \currentijn).&
\end{align*}
see Figure \ref{fig:caseiv}.
Also in this case the conclusion follows.

In the end we define 
\begin{equation} 
\elementXclip:=(\Gammaalp,\Gammabet,\Gammagam),\qquad \pGammaalpi:=\Gamma_{ij}^l=\Gamma_{ki}^r,\qquad i=1,2,3.
\end{equation}
The proof is achieved.

\end{proof}

From compactness of the space of BV connection, combining with Proposition \ref{Lem:AGcontinuous}, we see that the infimum in \eqref{Gammalimsupcurves} is attained. As a consequence, we can conclude the proof of   Theorem \ref{teo:main}.

\begin{Corollary}\label{COR}
We have
\begin{equation}\label{eqnminXclip}
\A(u,\sourcedisk)\leq |\sourcedisk| + \inf \{\AG(\elementXclip): \elementXclip \in \Xclip\}=|\sourcedisk| +\min \{\AG(\elementXclip): \elementXclip \in \Xc\}.
\end{equation}
\end{Corollary}

\section*{Acknowledgements}

 The present paper  benefits from the support of the GNAMPA (Gruppo Nazionale per l'Analisi Matematica, 
la Probabilit\`a e le loro Applicazioni) of INdAM 
(Istituto Nazionale di Alta Matematica).

\end{document}